\documentclass[preprint,10pt,authoryear]{elsarticle}

\usepackage[T1]{fontenc}																
\usepackage[utf8]{inputenc}															
\usepackage[english]{babel}															
\usepackage{amscd}																	
\usepackage{amsfonts}																
\usepackage{amsmath}																

\usepackage{amssymb}																

\usepackage{amsthm}																
\usepackage{bm}																	
\usepackage{bbm}																	
\usepackage{braket}																	
\usepackage[font=small,format=hang,labelfont={sf,bf}]{caption}									
\usepackage{changepage}															
\usepackage{comment}																
\usepackage[babel]{csquotes}															
\usepackage{empheq}																
\usepackage{enumitem}																
\usepackage{fancyhdr}																
\usepackage[flushmargin,hang,symbol]{footmisc} 
\usepackage{geometry}																
\usepackage{graphicx}																
\usepackage{guit}																	
\usepackage{indentfirst}																
\usepackage{lipsum}																	
\usepackage{listings}																
\usepackage{mathtools}																
\usepackage{mathrsfs}																
\usepackage{microtype}																
\usepackage{syntonly}
\usepackage{url}																	
\usepackage[english]{varioref}															
\usepackage{xcolor}																	
\usepackage{hyperref}																
	\hypersetup{
	colorlinks,
	citecolor=black,
	filecolor=black,
	linkcolor=black,
	urlcolor=black
	}

\newcommand{\N}{\mathbb{N}} 

\newcommand{\R}{\mathbb{R}}

\newcommand{\F}{\mathbb{F}}

\newcommand{\BB}{\mathcal{B}}
\newcommand{\CC}{\mathcal{C}}

\newcommand{\FF}{\mathcal{F}}

\newcommand{\HH}{\mathcal{H}}
\newcommand{\II}{\mathcal{I}}

\newcommand{\LL}{\mathcal{L}}

\newcommand{\NN}{\mathcal{N}}

\newcommand{\QQ}{\mathcal{Q}}

\newcommand{\UUU}{\mathscr{U}}

\renewcommand{\epsilon}{\varepsilon} 
\renewcommand{\theta}{\vartheta}
\renewcommand{\rho}{\varrho}
\renewcommand{\phi}{\varphi}

\newcommand{\T}{\mathsf{T}} 

\renewcommand{\u}{\bm{u}} 
\newcommand{\Py}{\bm{\Pi}} 

\newcommand{\D}{\textup{D}} 
\renewcommand{\O}{\textup{O}} 

\newcommand{\1}{\mathbbm{1}} 

\newcommand{\eqnot}{\coloneqq} 

\DeclarePairedDelimiter{\rounds}{{\mathopen{(}}}{{\mathclose{)}}} 
\DeclarePairedDelimiter{\bracks}{\lbrack}{\rbrack} 
\DeclarePairedDelimiter{\braces}{\lbrace}{\rbrace} 
\DeclarePairedDelimiter{\abs}{\lvert}{\rvert} 
\DeclarePairedDelimiter{\norma}{\lVert}{\rVert} 

\newcommand{\Lped}[3]{L_{#3}^{#1}\rounds{#2}} 
\newcommand{\LLped}[3]{\LL_{#3}^{#1}\rounds{#2}} 

\DeclareMathOperator{\probP}{\mathbf{P}} 
\DeclareMathOperator{\E}{\mathbf{E}} 
\DeclareMathOperator{\VaR}{VaR} 

\DeclareMathOperator{\meas}{\mathbf{m}} 

\DeclareMathOperator*{\argmin}{arg\,min} 

\newenvironment{sistema}
{\left\lbrace\begin{array}{@{}l@{}}}
{\end{array}\right.} 

\theoremstyle{plain} 
	\newtheorem{theorem}{Theorem} 
	\newtheorem{lemma}{Lemma} 
	
	\newtheorem{corollary}{Corollary} 
	
	\newtheorem{proposition}{Proposition} 

\theoremstyle{definition} 
	\newtheorem{definition}{Definition} 
	
	\newtheorem{assumption}{Assumption} 
	
	\newtheorem{example}{Example}
	
	\newtheorem{notation}{Notation}
	
	\newtheorem{problem}{Problem} 

\theoremstyle{remark} 
	\newtheorem{remark}{Remark} 
	
\journal{Journal of Mathematical Analysis and Applications}

\begin{document}

\begin{frontmatter}







\title{Subgame\:\!-perfect equilibrium strategies \\ for time\:\!-inconsistent recursive stochastic control problems}


\author{Elisa Mastrogiacomo}
\ead{elisa.mastrogiacomo@uninsubria.it}
\ead[url]{https://www.uninsubria.it/hpp/elisa.mastrogiacomo}
            
\author{Marco Tarsia}
\ead{marco.tarsia@uninsubria.it}
\ead[url]{https://www.uninsubria.it/hpp/marco.tarsia}

\begin{abstract}
We study time\:\!-inconsistent recursive stochastic control problems, i.e., for which Bellman's principle of optimality does not hold. For this class of problems classical optimal controls may fail to exist, or to be relevant in practice, and dynamic programming is not easily applicable. Therefore, the notion of optimality is defined through a game\:\!-theoretic framework by means of subgame\:\!-perfect equilibrium: we interpret our preference changes which, realistically, are inconsistent over time, as players in a game for which we want to find a Nash equilibrium. The approach followed in our work relies on the stochastic (Pontryagin) maximum principle: we adapt the classical spike variation technique to obtain a characterization of equilibrium strategies in terms of a generalized second-order Hamiltonian function defined through pairs of backward stochastic differential equations, even in the multidimensional case. The theoretical results are applied in the financial field to finite horizon investment-consumption policies with non-exponential actualization. Here the existence of non-trivial equilibrium policies is also ascertained.
\end{abstract}




\begin{keyword}
Equilibrium strategy \sep maximum principle \sep portfolio management \sep recursive stochastic control problem \sep time\:\!-inconsistency.

\MSC[2020] 60H30 \sep 91A23 \sep 91B70 \sep 93E20.
\end{keyword}

\end{frontmatter}




\section*{Introduction}

In this paper, we study time\:\!-inconsistent recursive stochastic control problems where the notion of optimality is defined by means of subgame\:\!-perfect equilibrium. In recent years, there has been increasing attention to time-inconsistent control problems, especially for their applicability to mathematical finance and economics. Time\:\!-inconsistency means that the Bellman’s principle of optimality does not hold. In practise, a restriction of an optimal control for a specific initial pair on a later time interval might not be optimal for that corresponding initial pair. This happens, for instance, in dynamic utility maximization problems for investment-consumption strategies under non-exponential discounting. To handle that problem in a time\:\!-consisting way, we need to introduce a different concept of solutions instead of the classical one. One of the main approaches to dealing with time\:\!-inconsistency is to find equilibrium controls, which are within a game theoretic framework.

In a continuous-time setting, such controls have been rigorously introduced in~\cite{ekelandlazrak06} and~\cite{ekelandpirvu08}, later completed in~\cite{ekelandmbodjipirvu12}, and can be thought of as ``infinitesimally optimal via spike variation'': i.e., they are optimal with respect to a penalty represented by deviations during an infinitesimal amount of time.

In particular, in~\cite{ekelandpirvu08}, the authors apply the classical maximum principle theory of~\cite{yongzhou99} to deal with the linear Merton portfolio management problem in the context of pseudo\:\!-exponential actualization, introducing the concept of subgame\:\!-perfect equilibrium policy as a notion to ensure the time\:\!-consistency of the portfolio strategy---possibly not unique. They arrive at equivalent formulations in terms of PDEs, or even ODEs, and integral equations due to the special form of discounting.  
		
The first aim of this paper is to perform similar operations for a more general control problem in the context of recursive utilities and we apply our results to the financial sphere. The approach followed in our work is inspired by~\cite{ekelandpirvu08} and~\cite{hu17} and relies on the stochastic maximum principle; see also~\cite{peng90},~\cite{peng93} and~\cite{brianddelyonhupardouxstoica03}. We adapt the classical spike variation technique to obtain a characterization of closed-loop equilibrium strategies in terms of a generalized Hamiltonian function $\HH$ defined through a pair of backward stochastic differential equations (BSDEs). Our generalized Hamiltonian function, compared with the classical one, contains the driver coefficient of the recursive utility (which has more variables than its analogue in the classical case) and involves a second-order stochastic process.

We emphasize that, similarly to the classical case, equilibrium strategies are characterized through both a necessary condition and a sufficient condition involving the generalized Hamiltonian function; whereas, contrary to the classical case, this sufficient condition works even in the absence of extra convexity assumptions. Also, we point out that the spike variation technique, explicitly required by the definition of equilibrium policy, applies indiscriminately to the case in which the control domain $U$ satisfies particular geometric conditions such as convexity or linearity and to more general cases; see~\cite{jizhou06} for a different approach (a terminal perturbation method) which is applicable only to the case of time\:\!-consistent optimal controls in the classical sense. 

Later on, the theoretical results are applied in the financial field to finite horizon investment-consumption policies with non-exponential actualization (e.g., a hyperbolic one). In particular, we extend the results contained in the aforementioned works~\cite{ekelandlazrak06},~\cite{ekelandmbodjipirvu12} and~\cite{ekelandpirvu08} by introducing the recursive utilities. An explicit characterization will smoothly be computed for a portfolio problem that we will consider at the end. We would like to observe here that the assumptions we adopt for the coefficients of our stochastic control system---a decoupled forward\!\:-backward stochastic differential equation---are substantially those proposed in~\cite{hu17}, although the concept of equilibrium is not used there. 
We look for controls in feedback or closed-loop form, partially mimicking what is done in~\cite{ekelandpirvu08} (or~\cite{ekelandmbodjipirvu12}), where explicit calculations are feasible. The shape of the recursive utility follows the classic Uzawa type (see~\cite{duffieepstein92}), but many other choices are possible. See also~\cite{gundelweber08},~\cite{imkellerdosreis10},~\cite{yong12} and~\cite{pennerreveillac15}. 

The theory of recursive optimal control problems in continuous time has attracted remarkable attention in recent years, both from a theoretical and an applicative viewpoint. For the time\:\!-consistent framework, we refer in particular to the fundamental works~\cite{duffieepstein92} and~\cite{elkarouipengquenez01} (see also~\cite{elkarouipengquenez97} and references therein). Time\:\!-inconsistent problems were first analyzed through subgame\:\!-perfect equilibrium strategies by~\cite{strotz55} and~\cite{pollak68}, and this line of research has been pursued by many others. We mention the series of studies by Yong (see, e.g.,~\cite{weiyongyu17}, ), whose approach focuses on dynamic programming, i.e., on Hamilton--Jacobi--Bellman equations (HJB equations). In \cite{bjorkkhapkomurgoci17}, the time\:\!-inconsistent control problem is considered in a general Markov framework, and an extended HJB equation together with the verification theorem are derived. Considering the hyperbolic discounting,~\cite{ekelandmbodjipirvu12} studies the portfolio management problem for an investor who is allowed to consume and take out life insurance, and the equilibrium strategy is characterized by an integral equation. See also~\cite{hujinzhou12},~\cite{bjorkmurgoci14},~\cite{bjorkmurgocizhou14},~\cite{bjorkkhapkomurgoci17},~\cite{hujinzhou17},~\cite{yong12},~\cite{yong14} and references therein for various kinds of problems. 

Very recently,~\cite{HaSIAM21} considered a time\:\!-inconsistent investment-consumption problem with random endowments in a possibly incomplete market under general discount functions. 
Finally, we mention~\cite{Ha21MCRF} for recent results concerning with time\:\!-inconsistent recursive stochastic control problems where the cost functional is defined by the solution to a backward stochastic Volterra integral equation. Differently from our approach, anyway, the author here focuses on open-loop equilibrium controls. 

%

The contributions of this paper are summarized as follows.

\begin{itemize}
	
	\item[--] The extension to the framework of recursive stochastic control problems of the definition of the equilibrium concept (given in~\cite{ekelandpirvu08}) within a class of closed-loop strategies.

	\item[--] The characterization of an equilibrium policy through the solution of a flow of BSDEs and to show that, under sufficiently general assumptions of the coefficients, this flow of BSDEs has a solution.

	\item[--] The formulation of a necessary and sufficient condition for a closed-loop equilibrium control via variational methods (Theorem~\ref{th:maxprinc}).

	\item[--] The treatment of an illustrative example concerned with an investment-consumption problem under non-exponential discounting, which also shows the existence of an equilibrium policy.

\end{itemize}

The remainder of the paper is organized as follows. In Section~\ref{sec:notations}, we introduce the notation. In Section~\ref{sec:defassprob}, we formulate the notion of (subgame\:\!-perfect) equilibrium policy and the class of problems which we are interested in. In Section~\ref{sec:preliminary}, we recall some preliminary results. In Section~\ref{sec:unconstrained}, we present and prove necessary and sufficient conditions for the existence of an equilibrium policy, even with a brief mention to the multidimensional case, deferring the more technical part to~\ref{subsec:lemmaliminf}. In Section~\ref{sec:applicationuncon}, we analyze a significant portfolio management problem as an application of the results obtained in the previous sections. In Section~\ref{sec:conclusions}, we discuss possible future research.




\section{Notation}
\label{sec:notations}


Set $T\in\mathopen{]}0,\infty\mathclose{[}$ as a finite deterministic \emph{horizon} and let $\rounds{\Omega,\FF,\probP}$ be a complete probability space such that we can define a one\:\!-dimensional Brownian motion, or Wiener process, $W=(W(t))_{t\in\mathopen{[}0,T\mathclose{]}}$ on it. Let $\F=(\FF_t)_{t\in\mathopen{[}0,T\mathclose{]}}$ be the completed filtration generated by $W$, for which we suppose that
\[
\FF_T=\FF
\]
(system noise is the only source of uncertainty in the problem). Thus, the filtered space $\rounds{\Omega,\FF,\F,\probP}$ satisfies the usual conditions. In this regard, see, e.g.,~\cite[Chap. 1, Sect. 2]{yongzhou99}.

\begin{remark}
For any non\:\!-\:\!empty set $\II$ of indices $\imath$, we will keep implicit the dependence on the sample variable $\omega\in\Omega$ for each stochastic process on $\II\times\Omega$, as is usually done (and as indeed we have just done for $W$). We specify also that any stochastic process on $\II\times\Omega$ must be seen as its equivalence class given by the quotient with respect to the equivalence relation $\sim$ of indistinguishability: i.e., for any process $X=(X(\imath))_{\imath\in\II}$, $\tilde{X}=(\tilde{X}(\imath))_{\imath\in\II}$ on $\II\times\Omega$, $X \sim \tilde{X}$ if and only if
\[
\probP\bracks[\big]{\forall\;\imath\in\II\text{, }X(\imath) = \tilde{X}(\imath)} = 1.
\]
\end{remark}

We introduce the following, rather familiar, notation, in which $E\in\BB(\R)$, $d\in\N^*\:\!\!=\N\setminus\Set{\!0\!}$, $V\:\!\! \subseteq \R^d \:\!\! = \R^{d \times 1}$\:\!\! is a vector subspace, $t,\tau\in\mathopen{[}0,T\mathclose{]}$ with $t \neq T$, and $p\in{[}1,\infty\mathclose{[}$.

\begin{itemize}[leftmargin=*]

	\item $\lesssim$\:\!. Less than or equal to, unless there are positive multiplicative constants (independent of what is involved) about which we are not particularly interested in being more explicit.

	\item $\meas\bracks{\,\bm{\cdot}\,}$. The one\:\!-dimensional Lebesgue measure on $\BB(\R)$.
	
	\item $\1_{\:\!\!E}\rounds{\:\!\bm{\cdot}\:\!}$. The indicator function of the set $E$, i.e., for any $x\in\R$,
\[
\1_{\:\!\!E}(x) \eqnot
\begin{cases}
1, & \text{if $x \in E$,} \\
0, & \text{if $x \notin E$.}
\end{cases}
\]
	
	\item $\E\bracks{\,\bm{\cdot}\,}$. The expected value w.r.t. $\probP$ of a ($\probP$-\:\!integrable) real\:\!-valued $\FF$-\:\!measurable random variable $X$ on $\Omega$, i.e.,
\[
\E\bracks{X} \equiv \E X \eqnot \int_\Omega\:\!\! X(\omega)\,d\:\!\!\probP(\omega) \in \R.
\]
	
	\item $\Lped{p}{\Omega;V}{\tau}$. The (Banach) space of $V$-\:\!valued $\FF_\tau$-\:\!measurable random variables $X$ on $\Omega$ such that
\[
\norma{X}_p^p \eqnot \E{\abs{X}}^p\:\!\! < \infty.
\]
	
	\item $\Lped{\infty}{\Omega;V}{\tau}$. The space of ($\probP$-\:\!a.s.) bounded $V$-\:\!valued $\FF_\tau$-\:\!measurable random variables $X$ on $\Omega$, i.e., with
\[
\norma{X}_\infty \eqnot \inf\Set{\:\!\! K\in\mathopen{[}0,\infty\mathclose{[} \mid \abs{X} \leq K \text{ $\probP$-\:\!a.s.} \:\!\!} < \infty
\]
where eventually, by convention,
\[
\inf\emptyset = \infty.
\]
	
	\item $\LLped{p}{t,T;V}{\F}$. The space of $V$-\:\!valued $(\FF_s)_{s\in\mathopen{[}t,T\mathclose{]}}$-\:\!progressively measurable processes $X\:\!\!=\,\!\!(X(s))_{s\in\mathopen{[}t,T\mathclose{]}}$ on $\mathopen{[}t,T\mathclose{]}\times\Omega$ (or $X(\:\!\bm{\cdot}\:\!)$, for short) such that
\[
\norma{X(\:\!\bm{\cdot}\:\!)}_p^p \equiv \norma{X}_p^p \eqnot \E\int_t^T\!{\abs{X(s)}}^p\:\!ds < \infty.
\]
	
	\item $\LLped{\infty}{t,T;V}{\F}$. The space of bounded $V$-\:\!valued $(\FF_s)_{s\in\mathopen{[}t,T\mathclose{]}}$-\:\!progressively measurable processes $X=(X(s))_{s\in\mathopen{[}t,T\mathclose{]}}$ on $\mathopen{[}t,T\mathclose{]}\times\Omega$, i.e., with
\[
\norma{X(\:\!\bm{\cdot}\:\!)}_\infty \equiv \norma{X}_\infty \eqnot \inf\Set{\! K\in\mathopen{[}0,\infty\mathclose{[} | \textstyle{\sup_{s\in\mathopen{[}t,T\mathclose{]}} \abs{X(s)}} \leq K \text{ $\probP$-\:\!a.s.} \!} < \infty.
\]
	
	\item $\LLped{p}{\Omega;\CC\rounds{\mathopen{[}t,T\mathclose{]};V}}{\F}$. The space of $V$-\:\!valued $(\FF_s)_{s\in\mathopen{[}t,T\mathclose{]}}$-\:\!adapted ($\probP$-\:\!a.s.) continuous processes $X=(X(s))_{s\in\mathopen{[}t,T\mathclose{]}}$ on $\mathopen{[}t,T\mathclose{]}\times\Omega$ such that
\[
\norma{X(\:\!\bm{\cdot}\:\!)}_{\CC,\:\!p}^p \equiv \norma{X}_{\CC,\:\!p}^p \eqnot \E\:\!\!\sup_{s\in\mathopen{[}t,T\mathclose{]}}\:\!\!{\abs{X(s)}}^p\:\!\! < \infty.
\]

\end{itemize}

\begin{remark}
\label{rem:normap}
For $X=(X(s))_{s\in\mathopen{[}t,T\mathclose{]}} \in \LLped{p}{t,T;V}{\F}$,
\[
\E\:\!\!{\mathopen{\bigg{(}}\int_t^T\!\abs*{X(s)}\:\!ds\mathclose{\bigg{)}}}^{\!p} \leq \norma{X(\:\!\bm{\cdot}\:\!)}_p^p < \infty
\]
(by the classical Jensen's inequality, w.r.t. the Lebesgue measure on $\mathopen{[}t,T\mathclose{]}$, for a.a. fixed $\omega\in\Omega$).
\end{remark} 

\begin{remark}
A $\R^d$-valued $(\FF_s)_{s\in\mathopen{[}t,T\mathclose{]}}$-\:\!adapted process $X=(X(s))_{s\in\mathopen{[}t,T\mathclose{]}}$ on $\mathopen{[}t,T\mathclose{]}\times\Omega$ admits an $(\FF_s)_{s\in\mathopen{[}t,T\mathclose{]}}$-\:\!progressively measurable modification (stochastically equivalent process) and, if $X$ is ($\probP$-\:\!a.s.) left continuous or right continuous as a process, then $X$ itself is $(\FF_s)_{s\in\mathopen{[}t,T\mathclose{]}}$-\:\!progressively measurable. See, e.g.,~\cite[Chap. 1, Sect. 2]{yongzhou99}.
\end{remark}


\section{Problem formulation} 
\label{sec:defassprob}


Before establishing definitions, assumptions and purposes, we would like to emphasize that the generalized optimization problem which we will deal with (see Problem~\ref{prob:UC}) would generally be affected by time\:\!-inconsistency due to the form of what will be our recursive utility (see Definition~\ref{def:recutility}) and, consequently, of the corresponding utility functional (see Definition~\ref{def:utilityfun}). Indeed, from a financial point of view, we can interpret it as having a structure of non\:\!-\:\!exponential time discounting and, therefore, a not constant (psychological) discount rate. This, on the other hand, explains why (spike variation technique and) equilibrium strategies are considered (see Definitions~\ref{def:equilibrium} and~\ref{def:equil4tuple}). See Section~\ref{sec:applicationuncon} for more details. We also refer to~\cite{ekelandpirvu08} and~\cite{weiyongyu17}. 

Take $n \in \N^*$\:\!\! and $\R^n$\:\!\! equipped with the Euclidean topology and the Borel $\sigma$-\:\!algebra $\BB(\R^n)$ with its Lebesgue measure, as will be the case for any other Euclidean space, and choose a \emph{control domain}
\[
U \in \BB(\R^n) \setminus \Set{\!\emptyset\!}\!,
\]
not necessarily convex nor bounded (in $\R^n$\:\!\!).

\begin{definition}[Admissible control]
\label{def:admissible}
Fix an arbitrary \emph{initial instant}
\[
t\in\mathopen{[}0,T\mathclose{[}.
\]
For an appropriate
\[
p\in\mathopen{[}2,\infty\mathclose{[}\:\!,
\]
that we do not want to give in explicit form (see~\cite[Introduction]{hu17}), we set
\begin{equation}
\label{eq:admissible}
\UUU\bracks{t,T} \doteq \Set{\u(\:\!\bm{\cdot}\:\!)\in\LLped{p}{t,T;\R^n}{\F} \mid \text{$\u(\:\!\bm{\cdot}\:\!)$ \emph{is} $U$-\emph{valued}}}
\end{equation} 
and we call an \emph{admissible control} any element $\u(\:\!\bm{\cdot}\:\!)$ of $\UUU\bracks{t,T}$.
\end{definition} 

\begin{remark}
\label{rem:Ubounded}
If $U$ is bounded, then the class $\UUU\bracks{t,T}$ in~\eqref{eq:admissible} simply coincides with the one constituted by the $U$-valued processes $\u(\:\!\bm{\cdot}\:\!)$ on $\mathopen{[}t,T\mathclose{]}\times\Omega$ such that $\u(\:\!\bm{\cdot}\:\!)$ is $(\FF_s)_{s\in\mathopen{[}t,T\mathclose{]}}$-\:\!progressively measurable.
\end{remark} 

\begin{definition}[Spike variation]
\label{def:spikevar}
Fix $t\in\mathopen{[}0,T\mathclose{[}$, $\bar{\u}(\:\!\bm{\cdot}\:\!)\in\UUU\bracks{t,T}$, $\epsilon\in\mathopen{]}0,T-t\mathclose{[}$, and a set
\[
E_t^{\:\!\epsilon}\:\!\! \in \BB(\mathopen{[}t,T\mathclose{]})
\]
with length
\[
\abs{E_t^{\:\!\epsilon}} \eqnot \meas\bracks{E_t^{\:\!\epsilon}} = \epsilon.
\]
For $\u(\:\!\bm{\cdot}\:\!)\in\UUU\bracks{t,T}$, we call the \emph{spike} (or \emph{needle}) \emph{variation of $\bar{\u}(\:\!\bm{\cdot}\:\!)$ w.r.t. $\u(\:\!\bm{\cdot}\:\!)$ and $E_t^{\:\!\epsilon}$\:\!\!} the admissible control $\bar{\u}^\epsilon(\:\!\bm{\cdot}\:\!)\in\UUU\bracks{t,T}$ defined by setting
\begin{equation} 
\label{eq:spikevar}
\bar{\u}^\epsilon\:\!\! \doteq \bar{\u} + (\u - \bar{\u})\1_{\:\!\!E_t^{\:\!\epsilon}}.
\end{equation} 
\end{definition} 

\begin{remark}
The spike variation $\bar{\u}^\epsilon(\:\!\bm{\cdot}\:\!)$ in~\eqref{eq:spikevar} is explicitly given, for any $s\in\mathopen{[}t,T\mathclose{]}$, by
\[
\bar{\u}^\epsilon(s) =
\begin{cases}
\bar{\u}(s), & \text{if $s \in \mathopen{[}t,T\mathclose{]} \setminus E_t^{\:\!\epsilon}$\:\!\!,} \\
\u(s), & \text{if $s \in E_t^{\:\!\epsilon}$\:\!\!.}
\end{cases}
\]
\end{remark}

\begin{notation}
Any alphabetic letter appearing as a subscript of a prescribed set that appears explicitly as (part of) the domain of a function such as, among others, $\u=(u_1,\dots,u_n)^\T$\:\!\! for
\[
U=U_{\u}
\]
or $x$ for
\[
\R = \R_x 
\]
should be seen as our preferred notation for the generic variable element of that domain.
\end{notation}

Fix four deterministic maps
\[
b,\sigma \colon \mathopen{[}0,T\mathclose{]}{_s} \times \R_x \times U_{\:\!\!\u} \to \R,
\]
\[
f \colon \mathopen{[}0,T\mathclose{]}{_s} \times \R_x \times U_{\:\!\!\u} \times \R_y \times \R_z \times \mathopen{[}0,T\mathclose{[}{_t} \to \R,
\]
and
\[
h \colon \R_x \times \mathopen{[}0,T\mathclose{[}{_t} \to \R
\]
(i.e., they do not depend on $\omega \in \Omega$) such that the following assumption holds, similarly to~\cite{hu17} (to which we refer to better understand why such strong assumptions are needed).

\begin{assumption}
\label{ass:assumption}
The maps $b,\sigma,f,h$ are continuous w.r.t. all their variables and, for any $t\in\mathopen{[}0,T\mathclose{[}$, there exists $L_t\in\mathopen{]}0,\infty\mathclose{[}$ such that, whatever map
\[
\phi(s,x,\u,y,z;t)
\]
between $b(s,x,\u)$, $\sigma(s,x,\u)$, $f(s,x,\u,y,z;t)$, and $h(x;t)$ is taken, and for any $s\in\mathopen{[}t,T\mathclose{]}$, $\u \in U$ and $x,y,z\in\R$,
\[
\abs*{\phi(s,x,\u,y,z;t)} \leq L_t\mathopen{\big{(}}1+\abs{x}+\abs{\u}+\abs{y}+\abs{z}\mathclose{\big{)}}.
\]
Next, $b,\sigma,h$ are of (differentiability) class $C^2$\:\!\! w.r.t. the variable $x\in\R$; $b_x,b_{xx},\sigma_x,\sigma_{xx}$ are bounded (on $\mathopen{[}t,T\mathclose{]}{_s} \times \R_x \times U_{\:\!\!\u}$) and continuous w.r.t. $\rounds{x,\u} \in \R \times U$; $h_x,h_{xx}$ are bounded and continuous (on $\R_x$); $f(\:\!\bm{\cdot}\,;t)$ is of class $C^2$\:\!\! w.r.t. $\rounds{x,y,z}\in\R^3$\:\!\!, with $\D f(\:\!\bm{\cdot}\,;t)$ and $\D^2 f(\:\!\bm{\cdot}\,;t)$ (gradient and Hessian matrix of $f(\:\!\bm{\cdot}\,;t)$ w.r.t. $\rounds{x,y,z}$ respectively) being bounded (on $\mathopen{[}t,T\mathclose{]}{_s} \times \R_x \times U_{\:\!\!\u} \times \R_y \times \R_z$) and continuous w.r.t. $\rounds{x,\u,y,z}\in\R \times U \times \R \times \R$.
\end{assumption} 

\begin{remark}
Regarding Assumption~\ref{ass:assumption}, we point out the following.
\begin{itemize}[leftmargin=*]


	\item The relations with $\phi=b$ and $\phi=\sigma$ could really depend on $t$ through the fact that $s\in\mathopen{[}t,T\mathclose{]}$. 
	
	\item The conditions of sublinear growth and boundedness imply something that is somehow stronger than implied by the classic conditions in~\cite[Chap. 3, Sect. 3]{yongzhou99}: more precisely, for any $t\in\mathopen{[}0,T\mathclose{[}$, there exists $L_t\in\mathopen{]}0,\infty\mathclose{[}$ such that, whatever map $\phi(s,x,\u,y,z;t)$ between $b(s,x,\u)$, $\sigma(s,x,\u)$, $f(s,x,\u,y,z;t)$, and $h(x;t)$ is taken, for any $s\in\mathopen{[}t,T\mathclose{]}$, $\u,\hat{\u} \in U$, and $x,\hat{x},y,\hat{y},z,\hat{z}\in\R$,
\begin{multline*}
\abs*{\phi(s,x,\u,y,z;t)-\phi(s,\hat{x},\hat{\u},\hat{y},\hat{z};t)} \\[0.5ex]
	\lor \abs*{\phi_x(s,x,\u,y,z;t)-\phi_x(s,\hat{x},\hat{\u},\hat{y},\hat{z};t)} \qquad \qquad \qquad \quad \\[0.5ex]
		\lor \abs*{\phi_{xx}(s,x,\u,y,z;t)-\phi_{xx}(s,\hat{x},\hat{\u},\hat{y},\hat{z};t)} \\[0.5ex]
			\leq L_t\mathopen{\big{(}}\:\!\abs{x-\hat{x}}+\abs{\u-\hat{\u}}+\abs{y-\hat{y}}+\abs{z-\hat{z}}\:\!\mathclose{\big{)}}
\end{multline*}
and
\[
\abs*{\phi(s,0,\u,0,0;t)} \leq L_t. 
\]
	
	
	\item See Section~\ref{sec:applicationuncon} for a situation where the maps $b$, $\sigma$, $f$ and $h$ satisfy all the regularity conditions required by Assumption~\ref{ass:assumption}.
	
\end{itemize}
\end{remark}

Now, choose a \emph{state domain}
\[
I \subseteq \R,
\]
that is a non-empty open interval. For fixed $t\in\mathopen{[}0,T\mathclose{[}$, all the following stochastic differential equations and corresponding (adapted) solutions are taken on
\[
\mathopen{[}t,T\mathclose{]}\times\Omega.
\]

\begin{definition}[Admissible state process]
\label{def:state}
Fix $t\in\mathopen{[}0,T\mathclose{[}$ and $\u(\:\!\bm{\cdot}\:\!)\in\UUU\bracks{t,T}$. For an arbitrary fixed \emph{initial state}
\[
x \in I,
\]
we call the \emph{state equation} or \emph{controlled system} (in the \emph{strong formulation}) the forward stochastic differential equation
\begin{equation}
\label{eq:state}
\begin{sistema}
dX(s) = b(s,X(s),\u(s))\:\!ds + \sigma(s,X(s),\u(s))\:\!dW(s), \\[0.5ex]
X(t) = x,
\end{sistema}
\end{equation} 
(where $s\in\mathopen{[}t,T\mathclose{]}$) and we call an \emph{admissible state process} any solution $X(\:\!\bm{\cdot}\:\!)$ of~\eqref{eq:state} that belongs to
\[
\LLped{2}{t,T;\R}{\F}.
\]
\end{definition} 

\begin{remark}
Regarding Definition~\ref{def:state}, we point out the following.
\begin{itemize}[leftmargin=*]
	\item The equation~\eqref{eq:state} is a controlled forward stochastic differential equation (FSDE) in It\^o differential form, with finite deterministic horizon $T$ and with random coefficients that depend on the sample $\omega\in\Omega$ only through $\u(\:\!\bm{\cdot}\:\!)$ and $X(\:\!\bm{\cdot}\:\!)$ itself, and it depends also on $b$, $\sigma$, $t$, and $x$ (as well as $W$).
	
	\item The term ``strong formulation'', which henceforth we will not repeat, alludes to the fact that the filtered space of probability $\rounds{\Omega,\FF,\F,\probP}$ is fixed a priori together with $W$ and  therefore must not be sought as part of the solution of~\eqref{eq:state} (see, e.g.,~\cite[Chap. 1, Sect. 6]{yongzhou99}).
	
	\item Suppose Assumption~\ref{ass:assumption} holds, at least as regards $b$ and $\sigma$. Fix $t\in\mathopen{[}0,T\mathclose{[}$, $x \in I$, and $\u(\:\!\bm{\cdot}\:\!)\in\UUU\bracks{t,T}$. Then there exists a unique solution
	\[
	X(\:\!\bm{\cdot}\:\!) \in \LLped{2}{\Omega;\CC\rounds{\mathopen{[}t,T\mathclose{]};\R}}{\F}
	\]
	of the FSDE~\eqref{eq:state} and its It\^o integral form is given, for any $s\in\mathopen{[}t,T\mathclose{]}$, by
	\begin{equation}
	\label{eq:X}
	X(s) = x + \int_t^{\:\!s}\! b(r,X(r),\u(r))\,dr + \int_t^{\:\!s}\! \sigma(r,X(r),\u(r))\,dW(r).
	\end{equation} 
	See, e.g., Proposition~\ref{prop:fbsde} below.
	
\end{itemize}
\end{remark}

\begin{notation}
\label{not:X}
We specify the dependence on the elements involved above by writing
\[
X(\:\!\bm{\cdot}\:\!) = X^{\:\!t,x,\u}(\:\!\bm{\cdot}\:\!) \eqnot X(\:\!\bm{\cdot}\,;t,x,\u(\:\!\textbf{-}\:\!))
\]
(see Definition~\ref{def:state}).
\end{notation} 

\begin{assumption}
\label{ass:assumptionX}
For $t\in\mathopen{[}0,T\mathclose{[}$, $\u(\:\!\bm{\cdot}\:\!)\in\UUU\bracks{t,T}$, $x \in I$, $s\in\mathopen{[}t,T\mathclose{]}$, and $\probP$-\:\!a.s.,
\[
X^{\:\!t,x,\u}(s) \in I.
\]
\end{assumption} 

\begin{remark}
Under Assumption~\ref{ass:assumptionX}, the interval $I$ may depend on $T$: e.g., the larger $T$ is, the larger $I$ may also be (however, in the worst case scenario, we could always take $I = \R$). Moreover, if we prefer, we can imagine that the domain component in the variable $x$ of the maps $b$, $\sigma$, $f$ and $h$ is restricted precisely to $I$ in such a way that the (analogue of) Assumption~\ref{ass:assumption} still holds. See Section~\ref{sec:applicationuncon} for a situation where Assumption~\ref{ass:assumptionX} is satisfied.
\end{remark}

\begin{definition}[Admissible recursive utility]
\label{def:recutility}
For $t\in\mathopen{[}0,T\mathclose{[}$, $x \in I$ and $\u(\:\!\bm{\cdot}\:\!)\in\UUU\bracks{t,T}$, let $X(\:\!\bm{\cdot}\:\!)$ be an admissible state process as in Definition~\ref{def:state}. We call a \emph{recursive (dis)utility system} a backward stochastic differential equation as
\begin{equation}
\label{eq:recutility}
\begin{sistema}
dY(s\:\!;t) = - \:\! f(s,X(s),\u(s),Y(s\:\!;t),Z(s\:\!;t);t)\:\!ds + Z(s\:\!;t)\:\!dW(s), \\[0.5ex] 
Y(T\:\!;t) = h(X(T);t),
\end{sistema}
\end{equation} 
(where $s\in\mathopen{[}t,T\mathclose{]}$) and we call an \emph{admissible recursive utility} any process $Y(\:\!\bm{\cdot}\,;t)$ such that $\rounds{Y(\:\!\bm{\cdot}\,;t),Z(\:\!\bm{\cdot}\,;t)}$ is a pair solution of~\eqref{eq:recutility} that belongs to
\[
\LLped{2}{\Omega;\CC\rounds{\mathopen{[}t,T\mathclose{]};\R}}{\F} \times \LLped{2}{t,T;\R}{\F}.
\]
\end{definition} 

\begin{remark}
Regarding Definition~\ref{def:recutility}, we point out the following.
\begin{itemize}[leftmargin=*]
	\item The equation~\eqref{eq:recutility} is a backward stochastic differential equation (BSDE) in It\^o differential form, decoupled from the FSDE~\eqref{eq:state} of Definition~\ref{def:state} on which it totally depends. Here, we have something much more general than a stochastic differential (dis)utility (SDU) in its original meaning: that is, essentially, a BSDE such as
	\[
	\begin{sistema}
	d\;\!\Xi(s\:\!;t) = -\:\!F(s,\u(s),\Xi(s\:\!;t);t)\:\!ds + \O(s\:\!;t)\:\!dW(s), \\[0.5ex]
	\Xi(T\:\!;t) = \xi_{\:\!t},
	\end{sistema}
	\]
	(where $\xi_{\:\!t} \in \Lped{2}{\Omega;\R}{T}$ and $s\in\mathopen{[}t,T\mathclose{]}$). See, e.g.,~\cite{duffieepstein92} and references therein.
	
	\item The term ``disutility'' anticipates the fact that there will be something to be minimized (not maximized). This will be done in a more general sense than the classic one: precisely, in the sense of subgame\:\!-perfect equilibrium strategies.
	
	\item Suppose Assumption~\ref{ass:assumption} holds and fix $t\in\mathopen{[}0,T\mathclose{[}$, $x \in I$, and $\u(\:\!\bm{\cdot}\:\!)\in\UUU\bracks{t,T}$. Then there exists a unique pair solution
	\[
	\rounds{Y(\:\!\bm{\cdot}\,;t),Z(\:\!\bm{\cdot}\,;t)} \in \LLped{2}{\Omega;\CC\rounds{\mathopen{[}t,T\mathclose{]};\R}}{\F} \times \LLped{2}{t,T;\R}{\F}
	\]
	of the BSDE~\eqref{eq:recutility}, whose It\^o integral form is given, for any $s\in\mathopen{[}t,T\mathclose{]}$, by
	\begin{equation}
	\label{eq:YZ}
	Y(s\:\!;t) = h(X(T);t) + \int_{\:\!\!s}^T\! f(r,X(r),\u(r),Y(r\:\!;t),Z(r\:\!;t);t)\,dr - \int_{\:\!\!s}^T\! Z(r\:\!;t)\,dW(r).
	\end{equation} 
	See Proposition~\ref{prop:fbsde} below and, for everything related to the fundamental theory of BSDEs, see, e.g.,~\cite[Chap. 7]{yongzhou99}  and~\cite{elkarouipengquenez97}.
	
\end{itemize}
\end{remark}

\begin{notation}
\label{not:YZ}
We specify the dependence on the elements involved above by writing
\[
Y(\:\!\bm{\cdot}\,;t) = Y^{\:\!x,\u}(\:\!\bm{\cdot}\,;t) \eqnot Y(\:\!\bm{\cdot}\,;t,x,\u(\:\!\textbf{-}\:\!))
\]
and
\[
Z(\:\!\bm{\cdot}\,;t) = Z^{\:\!x,\u}(\:\!\bm{\cdot}\,;t) \eqnot Z(\:\!\bm{\cdot}\,;t,x,\u(\:\!\textbf{-}\:\!))
\]
(see Definition~\ref{def:recutility}).
\end{notation} 

\begin{remark}
\label{rem:deterministic}
$Y(t\:\!;t)$ is a deterministic constant. Indeed, since $x$ is a deterministic constant, $X(T)$ (see~\eqref{eq:X}) is measurable w.r.t. the completed $\sigma$\:\!-\:\!algebra $\tilde\FF_{t,T}$ on $\Omega$ generated by the process
\[
(W(s) - W(t))_{s\in\mathopen{[}t,T\mathclose{]}},
\]
and therefore $Y(t\:\!;t)$ (see~\eqref{eq:YZ}) is simultaneously measurable w.r.t. the two mutually independent $\sigma$\:\!-\:\!algebras $\FF_t$ and $\tilde\FF_{t,T}$ (an argument of this kind is found in, e.g.,~\cite{debusschefuhrmantessitore07}). Consequently,
\[
\begin{split}
\E Y(t\:\!;t) &\equiv \E\bracks*{\int_t^T\! f(s,X(s),\u(s),Y(s\:\!;t),Z(s\:\!;t);t)\,ds + h(X(T);t)} \\[1ex]
	&= Y(t\:\!;t).
\end{split}
\]
On the other hand, it is not possible to establish in general that $Z(t\:\!;t)$ is a deterministic constant.
\end{remark} 

\begin{remark}
We will prefer the notation $X(\:\!\bm{\cdot}\:\!)$ to other possible notations such as $X(\:\!\bm{\cdot}\,;t)$, but we will retain the notations $Y(\:\!\bm{\cdot}\,;t)$ and $Z(\:\!\bm{\cdot}\,;t)$.
\end{remark}

\begin{definition}[Recursive stochastic control problem]
\label{def:fbsde}
Fix $t\in\mathopen{[}0,T\mathclose{[}$, $x \in I$ and $\u(\:\!\bm{\cdot}\:\!)\in\UUU\bracks{t,T}$. We call the \emph{recursive stochastic control problem} the combination of the two stochastic differential equations~\eqref{eq:state} and~\eqref{eq:recutility}, i.e.,
\begin{equation}
\label{eq:fbsde}
\begin{sistema}
dX(s) = b(s,X(s),\u(s))\:\!ds + \sigma(s,X(s),\u(s))\:\!dW(s), \\[0.5ex]
dY(s\:\!;t) = - \:\! f(s,X(s),\u(s),Y(s\:\!;t),Z(s\:\!;t);t)\:\!ds + Z(s\:\!;t)\:\!dW(s), \\[0.5ex]
X(t) = x,\quad Y(T\:\!;t) = h(X(T);t),
\end{sistema}
\end{equation} 
(where $s\in\mathopen{[}t,T\mathclose{]}$) and, if $\rounds{X(\:\!\bm{\cdot}\:\!),Y(\:\!\bm{\cdot}\,;t),Z(\:\!\bm{\cdot}\,;t)}$ is a solution of~\eqref{eq:fbsde} that belongs to
\[
\LLped{2}{t,T;\R}{\F} \times \LLped{2}{\Omega;\CC\rounds{\mathopen{[}t,T\mathclose{]};\R}}{\F} \times \LLped{2}{t,T;\R}{\F},
\]
then we call $\rounds{\u(\:\!\bm{\cdot}\:\!),X(\:\!\bm{\cdot}\:\!),Y(\:\!\bm{\cdot}\,;t),Z(\:\!\bm{\cdot}\,;t)}$  an \emph{admissible 4-tuple}.
\end{definition} 

\begin{remark}
The equation/system~\eqref{eq:fbsde} is a controlled decoupled forward\:\!-backward stochastic differential equation/system (FBSDE) in It\^o differential form and, of course, we could use Notations~\ref{not:X} and~\ref{not:YZ} for the respective components of the corresponding solution.
\end{remark}

Regarding the recursive stochastic control problem~\eqref{eq:fbsde} of Definition~\ref{def:fbsde}, the following standard result for existence, uniqueness, and regularity holds (see, e.g.,~\cite{mamorelyong99}).
\begin{proposition}
\label{prop:fbsde}
Suppose Assumption~\ref{ass:assumption} holds and fix $t\in\mathopen{[}0,T\mathclose{[}$, $x \in I$ and $\u(\:\!\bm{\cdot}\:\!)\in\UUU\bracks{t,T}$. Then there exists a unique solution
\[
\rounds{X(\:\!\bm{\cdot}\:\!),Y(\:\!\bm{\cdot}\,;t),Z(\:\!\bm{\cdot}\,;t)} \in \LLped{2}{\Omega;\CC\rounds{\mathopen{[}t,T\mathclose{]};\R}}{\F}^2\:\!\! \times \LLped{2}{t,T;\R}{\F}
\]
of the FBSDE~\eqref{eq:fbsde} and
\[
\E\:\!\!\sup_{s\in\mathopen{[}t,T\mathclose{]}}\:\!\!{\abs{X(s)}}^2\:\!\! + \E\:\!\!\sup_{s\in\mathopen{[}t,T\mathclose{]}}\:\!\!{\abs{Y(s\:\!;t)}}^2\:\!\! + \E\int_t^T\!{\abs{Z(s\:\!;t)}}^2\:\!ds \:\! \lesssim \:\! 1 + x^2\:\!\! + \E\int_t^T\!{\abs{\u(s)}}^2\:\!ds.
\]
\end{proposition} 

\begin{definition}[Utility functional]
\label{def:utilityfun}
Suppose Assumption~\ref{ass:assumption} holds, and fix $t\in\mathopen{[}0,T\mathclose{[}$ and $x \in I$. For any $\u(\:\!\bm{\cdot}\:\!)\in\UUU\bracks{t,T}$, consider the solution $\rounds{X(\:\!\bm{\cdot}\:\!),Y(\:\!\bm{\cdot}\,;t),Z(\:\!\bm{\cdot}\,;t)}$ of the FBSDE~\eqref{eq:fbsde} as in Proposition~\ref{prop:fbsde}. We call \emph{(dis)utility} or \emph{cost functional} the map $J(\:\!\bm{\cdot}\,;t,x)\colon\UUU\bracks{t,T}\to\R$ given by
\begin{equation}
\label{eq:utilityfun}
J(\u(\:\!\bm{\cdot}\:\!);t,x) \doteq Y(t\:\!;t)
\end{equation} 
(see also Remark~\ref{rem:deterministic}).
\end{definition} 

\begin{remark}
The functional $J(\:\!\bm{\cdot}\,;t,x)$ in~\eqref{eq:utilityfun} is a real\:\!-valued generalized Bolza\:\!-type functional and, for any $\u(\:\!\bm{\cdot}\:\!)\in\UUU\bracks{t,T}$, it can be written as
\[
J(\u(\:\!\bm{\cdot}\:\!);t,x) = \E\bracks*{\int_t^T\! f(s,X(s),\u(s),Y(s\:\!;t),Z(s\:\!;t);t)\,ds + h(X(T);t)} \:\!\! .
\]
Note that the \emph{running} or \emph{intertemporal} \emph{utility} and the \emph{terminal utility} are explicitly specified in the above expression (see again Remark~\ref{rem:deterministic}). In particular, the more constant $f$ is w.r.t. the variables $y$ and $z$, the more we return to the classical sphere of stochastic optimal control theory.
\end{remark}

\begin{definition}[Equilibrium policy]
\label{def:equilibrium}
Suppose Assumptions~\ref{ass:assumption} and~\ref{ass:assumptionX} hold. We call a \emph{(subgame\:\!-perfect) equilibrium policy associated with $T,I,U,W$\:\!\! and $b,\sigma,f,h$} any measurable map
\[
\Py \colon \mathopen{[}0,T\mathclose{]}{_s} \times I_x \to U
\]
such that, for any $t\in\mathopen{[}0,T\mathclose{[}$ and $x \in I$, there exists a unique $I$\:\!-\:\!valued It\^o process
\begin{equation}
\label{not:XPy}
\bar{X}(\:\!\bm{\cdot}\:\!) \eqnot X^{t,x,\Py}(\:\!\bm{\cdot}\:\!)
\end{equation} 
that is a solution belonging to $\LLped{2}{t,T;\R}{\F}$ of the FSDE
\[
\begin{sistema}
dX(s) = b\bigl{(}s,X(s),\Py(s,X(s))\bigr{)}\:\!ds + \sigma\bigl{(}s,X(s),\Py(s,X(s))\bigr{)}\:\!dW(s), \\[0.5ex]
X(t) = x,
\end{sistema}
\]
(where $s\in\mathopen{[}t,T\mathclose{]}$) and is such that if we denote, for any $s\in\mathopen{[}t,T\mathclose{]}$ (and $\probP$-\:\!a.s.),
\begin{equation} 
\label{eq:equilibriumcontrol}
\bar{\u}(s) \doteq \Py(s,\bar{X}(s)),
\end{equation} 
then we have $\bar{\u}(\:\!\bm{\cdot}\:\!)\in\UUU\bracks{t,T}$ and, for any other $\u(\:\!\bm{\cdot}\:\!)\in\UUU\bracks{t,T}$,
\begin{equation} 
\label{eq:liminf}
\liminf_{\epsilon\downarrow0}\:\!\frac{J(\bar{\u}^\epsilon(\:\!\bm{\cdot}\:\!);t,x) - J(\bar{\u}(\:\!\bm{\cdot}\:\!);t,x)}{\epsilon} \geq 0,
\end{equation} 
where, for $\epsilon\downarrow0$, $\bar{\u}^\epsilon(\:\!\bm{\cdot}\:\!)$ is the spike variation of $\bar{\u}(\:\!\bm{\cdot}\:\!)$ w.r.t. $\u(\:\!\bm{\cdot}\:\!)$ and $E_t^{\:\!\epsilon}$\:\!\! is given by
\begin{equation}
\label{eq:Eteps}
E_t^{\:\!\epsilon}\:\!\! \eqnot \mathopen{[}t,t+\epsilon\mathclose{]}
\end{equation} 
(see also Definition~\ref{def:spikevar}).
\end{definition} 

\begin{remark}
The $\liminf_{\epsilon\downarrow0}$ in~\eqref{eq:liminf} will turn out to be an actual limit (see Lemma~\ref{lemma:liminf} in Section~\ref{sec:unconstrained}). 
\end{remark}

\begin{notation}
\label{not:XYZ}
With respect to the notations of Definition~\ref{def:equilibrium}, 
\[
\bar{X}(\:\!\bm{\cdot}\:\!) \eqnot X^{\:\!t,x,\bar{\u}}(\:\!\bm{\cdot}\:\!), \quad \bar{Y}(\:\!\bm{\cdot}\,;t) \eqnot Y^{\:\!x,\bar{\u}}(\:\!\bm{\cdot}\,;t), \quad \bar{Z}(\:\!\bm{\cdot}\,;t) \eqnot Z^{\:\!x,\bar{\u}}(\:\!\bm{\cdot}\,;t)
\]
and
\[
X^\epsilon(\:\!\bm{\cdot}\:\!) \eqnot X^{\:\!t,x,\bar{\u}^\epsilon}\!(\:\!\bm{\cdot}\:\!), \quad Y^\epsilon(\:\!\bm{\cdot}\,;t) \eqnot Y^{\:\!x,\bar{\u}^\epsilon}\!(\:\!\bm{\cdot}\,;t), \quad Z^\epsilon(\:\!\bm{\cdot}\,;t) \eqnot Z^{\:\!x,\bar{\u}^\epsilon}\!(\:\!\bm{\cdot}\,;t)
\]
(similarly to Notation~\ref{not:X} and Notation~\ref{not:YZ}).
\end{notation} 

\begin{definition}[Equilibrium control/pair/4\:\!-tuple]
\label{def:equil4tuple}
With respect to the notations of Definition~\ref{def:equilibrium} and Notation~\ref{not:XYZ}, we call:
\[
\bar{\u}(\:\!\bm{\cdot}\:\!)
\]
a \emph{(subgame\:\!-perfect) equilibrium control} (or \emph{strategy}),
\[
\rounds{\bar{\u}(\:\!\bm{\cdot}\:\!),\bar{X}(\:\!\bm{\cdot}\:\!)}
\]
a \emph{(subgame\:\!-perfect) equilibrium pair} and
\[
\rounds{\bar{\u}(\:\!\bm{\cdot}\:\!),\bar{X}(\:\!\bm{\cdot}\:\!),\bar{Y}(\:\!\bm{\cdot}\,;t),\bar{Z}(\:\!\bm{\cdot}\,;t)}
\]
a \emph{(subgame\:\!-perfect) equilibrium 4-tuple}.
\end{definition} 

\begin{remark}
\label{rem:liminfY}
Regarding Definitions~\ref{def:equilibrium} and~\ref{def:equil4tuple}, we point out the following.
\begin{itemize}[leftmargin=*]

	\item If $\bar{\u}(\:\!\bm{\cdot}\:\!)$ as in~\eqref{eq:equilibriumcontrol} is an optimal control in the classical sense, i.e., $\bar{\u}(\:\!\bm{\cdot}\:\!)$ minimizes the objective functional $J(\:\!\bm{\cdot}\,;t,x)$ over $\UUU\bracks{t,T}$, then $\bar{\u}(\:\!\bm{\cdot}\:\!)$ is also an equilibrium control (the opposite cannot be true, in general).

	\item We could specify an equilibrium policy/control/pair/4\:\!-tuple to be \emph{strong} in cases in which the inequality in~\eqref{eq:liminf} is strong, i.e., narrow (adjusting the entire sequel accordingly).

	\item  In general, we cannot expect an equilibrium policy/control/pair/4\:\!-tuple to be unique, even if it exists. It might therefore be an idea to select one uniquely through a \emph{constraint}: this will be the subject of our future studies (see also Section~\ref{sec:conclusions}).

	\item If $\rounds{\bar{\u}(\:\!\bm{\cdot}\:\!),\bar{X}(\:\!\bm{\cdot}\:\!)}$ is an equilibrium pair, then, by~\eqref{eq:equilibriumcontrol},
	\[
	\bar{\u}(t) = \Py(t,x)
	\]
	which, by definition of $\Py$, is a deterministic constant (vector in $U$).
	
	\item The condition~\eqref{eq:liminf} can be rewritten as
	\begin{equation}
	\label{eq:liminfY}
	\liminf_{\epsilon\downarrow0}\:\!\frac{Y^\epsilon(t\:\!;t) - \bar{Y}(t\:\!;t)}{\epsilon} \geq 0
	\end{equation} 
	(see Definition~\ref{def:utilityfun} and Notation~\ref{not:XYZ}).
	
	\item Section~\ref{sec:applicationuncon} attests that an equilibrium policy exists (at least in simple special cases).

\end{itemize}
\end{remark} 

The stochastic control problem that we will deal with, for which it is essentially a matter of obtaining a (Pontryagin) \emph{maximum principle}, can be stated as follows (see Definition~\ref{def:equil4tuple}).
\begin{problem}
\label{prob:UC}
Find \emph{necessary and sufficient conditions} for an equilibrium 4\:\!-tuple.
\end{problem} 


\section{Some preliminary results}
\label{sec:preliminary}


We recall a standard estimate for BSDEs that is decisive in~\cite{hu17} (on which we rely) and that can be found in, e.g.,~\cite{brianddelyonhupardouxstoica03}. Its meaning is essentially the continuous dependence of the pair solution on the assigned data.

\begin{lemma}
\label{lemma:continuitybsdes}
Fix $t\in\mathopen{[}0,T\mathclose{[}$. Take $p\in\mathopen{]}1,\infty\mathclose{[}$, $\xi_{\:\!t},\hat{\xi}_{\:\!t}\in\Lped{p}{\Omega;\R}{T}$ and two measurable maps
\[
F,\hat{F}\colon \mathopen{[}0,T\mathclose{]}{_s} \times \R_y \times \R_z \times \Omega_\omega \times \mathopen{[}0,T\mathclose{[}{_t} \to \R
\]
that are $s$\:\!-\:\!$\probP$-\:\!uniformly Lipschitz w.r.t. $\rounds{y,z}$ such that, for any $y,z\in\R$, $F(\:\!\bm{\cdot}\:\!,y,z;t)$ and $\hat{F}(\:\!\bm{\cdot}\:\!,y,z;t)$ are real $(\FF_s)_{s\in\mathopen{[}t,T\mathclose{]}}$-\:\!progressively measurable processes with $F(\:\!\bm{\cdot}\:\!,0,0;t),\hat{F}(\:\!\bm{\cdot}\:\!,0,0;t)\in\LLped{p}{t,T;\R}{\F}$. Consider the BSDEs with parameters $\rounds{-F,\xi_{\:\!t}}$ and $\rounds{-\hat{F},\hat{\xi}_{\:\!t}}$ respectively, i.e.,
\[
\begin{sistema}
d\;\!\Xi(s\:\!;t) = -\:\!F(s,\Xi(s\:\!;t),\O(s\:\!;t);t)\:\!ds + \O(s\:\!;t)\:\!dW(s), \\[0.5ex]
\Xi(T\:\!;t) = \xi_{\:\!t} \:\! ,
\end{sistema}
\]
and
\[
\begin{sistema}
d\;\!\hat{\Xi}(s\:\!;t) = -\:\!\hat{F}(s,\hat{\Xi}(s\:\!;t),\hat{\O}(s\:\!;t);t)\:\!ds + \hat{\O}(s\:\!;t)\:\!dW(s), \\[0.5ex]
\hat{\Xi}(T\:\!;t) = \hat{\xi}_{\:\!t} \:\! .
\end{sistema}
\]
(where $s\in\mathopen{[}t,T\mathclose{]}$). Then there exists a constant $K_p\in\mathopen{]}0,\infty\mathclose{[}$ such that
\begin{multline*}
\E\:\!\!\sup_{s\in\mathopen{[}t,T\mathclose{]}}\:\!\!{\abs*{\:\! \Xi(s\:\!;t) - \hat{\Xi}(s\:\!;t)}}^p\:\!\! + \E\:\!\!{\mathopen{\bigg{(}}\int_t^T\!{\abs*{\O(s\:\!;t) - \hat{\O}(s\:\!;t)}}^2\:\!ds\mathclose{\bigg{)}}}^{\!p/2}\!\! \:\! \leq \:\! K_p \E{\abs*{\xi_{\:\!t} - \hat{\xi}_{\:\!t}}}^p\:\!\! \\[0.5ex]
		+ K_p \E\:\!\!{\mathopen{\bigg{(}}\int_t^T\:\!\!\abs*{F(s,\Xi(s\:\!;t),\O(s\:\!;t);t) - \hat{F}(s,\Xi(s\:\!;t),\O(s\:\!;t);t)}\:\!ds\mathclose{\bigg{)}}}^{\!p\:\!}.
\end{multline*}
\end{lemma} 

\begin{remark}
Regarding Lemma~\ref{lemma:continuitybsdes}, we point out the following.
\begin{itemize}[leftmargin=*]

	\item The constant $K_p$ depends also on $t$, $T$, and the Lipschitz constants, but neither on $\xi_{\:\!t}(\:\!\bm{\cdot}\:\!), \hat{\xi}_{\:\!t}(\:\!\bm{\cdot}\:\!)$, $\Xi(\:\!\bm{\cdot}\,;t), \O(\:\!\bm{\cdot}\,;t)$ nor on $\hat{\Xi}(\:\!\bm{\cdot}\,;t), \hat{\O}(\:\!\bm{\cdot}\,;t)$.
	
	\item The main result underlying the whole theory of BSDEs is the classic representation theorem of integrable square continuous martingales, and thus it is crucial that the reference filtration remains the completed filtration $\F$ generated by $W$.
	
	
\end{itemize}
\end{remark}

We conclude the current section with a brief discussion of the classic comparison theorem for BSDEs, which, in the context of linearity, boils down to a simple observation, Remark~\ref{rem:comparison} below, which will be important for our discussion, especially because what we will call \emph{adjoint equations} will be linear BSDEs. See, e.g.,~\cite{elkarouipengquenez97}. 

\begin{proposition}
\label{prop:linearBSDEs}
Fix $t\in\mathopen{[}0,T\mathclose{[}$. Take
\[
\beta(\:\!\bm{\cdot}\,;t),\gamma(\:\!\bm{\cdot}\,;t) \in \LLped{\infty}{t,T;\R}{\F}
\]
and $\eta(\:\!\bm{\cdot}\,;t)$ such that
\[
\begin{sistema}
d\:\!\eta(s\:\!;t) = \eta(s\:\!;t)\bracks[\big]{\beta(s\:\!;t)\:\!ds + \gamma(s\:\!;t)\:\!dW(s)}, \\[0.5ex]
\eta(t\:\!;t) = 1,
\end{sistema}
\]
(where $s\in\mathopen{[}t,T\mathclose{]}$), i.e., explicitly,
\[
\eta(s\:\!;t) = \exp\braces[\bigg]{\int_t^{\:\!s}\:\!\! \bracks*{\beta(r\:\!;t) - \frac{\gamma^2(r\:\!;t)}{2}}\:\!dr + \int_t^{\:\!s}\! \gamma(r\:\!;t)\,dW(r)}.
\]
Then, for any $\alpha(\:\!\bm{\cdot}\,;t) \in \LLped{2}{t,T;\R}{\F}$ and $\xi_{\:\!t} \in \Lped{2}{\Omega;\R}{T}$, there exists a unique pair solution
\[
\rounds{\Xi(\:\!\bm{\cdot}\,;t),\O(\:\!\bm{\cdot}\,;t)} \in \LLped{2}{\Omega;\CC\rounds{\mathopen{[}t,T\mathclose{]};\R}}{\F} \times \LLped{2}{t,T;\R}{\F}
\]
of the BSDE
\begin{equation}
\label{eq:linearBSDEs1}
\begin{sistema}
d\;\!\Xi(s\:\!;t) = -\:\!\bracks[\big]{\alpha(s\:\!;t) + \beta(s\:\!;t)\:\!\Xi(s\:\!;t) + \gamma(s\:\!;t)\O(s\:\!;t)}\:\!ds + \O(s\:\!;t)\:\!dW(s), \\[0.5ex]
\Xi(T\:\!;t) = \xi_{\:\!t},
\end{sistema}
\end{equation} 
(where $s\in\mathopen{[}t,T\mathclose{]}$) and $\Xi(\:\!\bm{\cdot}\,;t)$ is the conditional expectation given, for any $s\in\mathopen{[}t,T\mathclose{]}$ (and $\probP$-\:\!a.s.), by
\begin{equation}
\label{eq:linearBSDEs2}
\Xi(s\:\!;t) = \eta^{-1}(s\:\!;t) \E\bracks*{\eta(T\:\!;t)\:\!\xi_{\:\!t} + \:\!\!\int_{\:\!\!s}^T\! \eta(r\:\!;t)\:\!\alpha(r\:\!;t)\,dr \,\Bigg{|}\, \FF_s}.
\end{equation} 
\end{proposition} 

\begin{remark}
\label{rem:comparison}
Regarding Proposition~\ref{prop:linearBSDEs}, we point out the following.
\begin{itemize}[leftmargin=*]

	\item  Since $\eta(\:\!\bm{\cdot}\,;t)>0$, it follows from~\eqref{eq:linearBSDEs2} that
	\begin{equation}
	\label{eq:comparison}
	\begin{sistema}
	\xi_{\:\!t} \geq 0 \\[0.5ex]
	\alpha(\:\!\bm{\cdot}\,;t) \geq 0
	\end{sistema}
	\qquad \Longrightarrow \qquad
	\Xi(\:\!\bm{\cdot}\,;t) \geq 0
	\end{equation} 
	(similarly with $\leq$ everywhere) and, moreover, the narrow inequality for $\Xi(\:\!\bm{\cdot}\,;t)$ holds even if only one of the two other inequalities is narrow: e.g., 
	\[
	\begin{sistema}
	\xi_{\:\!t} > 0 \\[0.5ex]
	\alpha(\:\!\bm{\cdot}\,;t) \geq 0
	\end{sistema}
	\qquad \Longrightarrow \qquad
	\Xi(\:\!\bm{\cdot}\,;t) > 0 \:\! .
	\]
	
	\item In general, for $\tau\in\mathopen{[}s,T\mathclose{]}$,
	\[
	\eta(\tau\:\!;t)\;\!\eta^{-1}(s\:\!;t) \neq \eta(\tau\:\!;s)
	\]
	(owing to the dependence on $t$ of $\beta(\:\!\bm{\cdot}\,;t)$ and $\gamma(\:\!\bm{\cdot}\,;t)$), and  therefore we can expect that, as processes,
	\[
	\Xi(s\:\!;t) \neq \E\bracks*{\eta(T\:\!;s)\:\!\xi_{\:\!t} + \:\!\!\int_{\:\!\!s}^T\! \eta(r\:\!;s)\:\!\alpha(r\:\!;t)\,dr \,\Bigg{|}\, \FF_s}
	\]
	(in which the latter term  differs from $\Xi(s\:\!;s)$ through the dependence on $t$ of $\alpha(\:\!\bm{\cdot}\,;t)$ and $\xi_{\:\!t}$).

\end{itemize}
\end{remark} 


\section{A maximum principle: necessary and sufficient conditions}
\label{sec:unconstrained}


In this section, we solve Problem~\ref{prob:UC} by adapting the calculations of~\cite{hu17} appropriately. Therefore, heuristics concerning the form of the adjoint equations/processes and their respective generalized Hamiltonian functions are not provided.

Let us now highlight the key difference w.r.t. the classical maximum principle. As will be shown in~\ref{subsec:lemmaliminf}, the utility functional $J(\:\!\bm{\cdot}\,;t,x)$ must be optimized in the (``weak'') sense of equilibrium policies: in particular, through the usual spike variation technique. So, on the one hand, we proceed along the lines developed in~\cite[Chap. 3, Sect. 4]{yongzhou99}. On the other hand, since $J(\:\!\bm{\cdot}\,;t,x)$ has a precise (``strong'') structure that derives from a recursive utility system and, thus, a BSDE, we extend the powerful techniques proposed in~\cite{hu17} to the time\:\!-inconsistent framework.

We suppose Assumptions~\ref{ass:assumption} and~\ref{ass:assumptionX} hold and we fix ($t\in\mathopen{[}0,T\mathclose{[}$, $x \in I$ and) an admissible 4\:\!-tuple
\[
\rounds{\bar{\u}(\:\!\bm{\cdot}\:\!),\bar{X}(\:\!\bm{\cdot}\:\!),\bar{Y}(\:\!\bm{\cdot}\,;t),\bar{Z}(\:\!\bm{\cdot}\,;t)}
\]
(see Definition~\ref{def:fbsde}) that we see as a candidate equilibrium 4\:\!-tuple (see Definitions~\ref{def:equilibrium} and~\ref{def:equil4tuple}).

\begin{notation}
\label{not:phi}
For any map $\phi(s,x,\u)$ between $b(s,x,\u)$, $\sigma(s,x,\u)$ and their derivatives up to second order, and for $s\in\mathopen{[}t,T\mathclose{]}$ (and $\probP$-\:\!a.s.), we write
\[
\phi(s) \eqnot \phi(s,\bar{X}(s),\bar{\u}(s))
\]
and, for $\u(\:\!\bm{\cdot}\:\!)\in\UUU\bracks{t,T}$,
\[
\delta\phi(s) \eqnot \phi(s,\bar{X}(s),\u(s)) - \phi(s),
\]
while similarly, for any map $\phi(s,x,\u,y,z;t)$ between $f(s,x,\u,y,z;t)$ and its derivatives up to second order, and for $s\in\mathopen{[}t,T\mathclose{]}$ (and $\probP$-\:\!a.s.), we write
\[
\phi(s\:\!;t) \eqnot \phi(s,\bar{X}(s),\bar{\u}(s),\bar{Y}(s\:\!;t),\bar{Z}(s\:\!;t);t).
\]
\end{notation} 


Regarding the following, let us keep in mind Notation~\ref{not:phi}.

\begin{definition}[$\kappa(\:\!\bm{\cdot}\,;t)$]
\label{def:kappa}
Associated with $\rounds{\bar{\u}(\:\!\bm{\cdot}\:\!),\bar{X}(\:\!\bm{\cdot}\:\!),\bar{Y}(\:\!\bm{\cdot}\,;t),\bar{Z}(\:\!\bm{\cdot}\,;t)}$, we define the process $\kappa(\:\!\bm{\cdot}\,;t)$ as the solution of the linear FSDE
\begin{equation}
\label{eq:kappa}
\begin{sistema}
d\:\!\kappa(s\:\!;t) = \kappa(s\:\!;t)\bracks[\big]{f_y(s\:\!;t)\:\!ds + f_z(s\:\!;t)\:\!dW(s)}, \\[0.5ex]
\kappa(t\:\!;t) = 1,
\end{sistema}
\end{equation} 
(where $s\in\mathopen{[}t,T\mathclose{]}$).
\end{definition} 

\begin{remark}
The process $\kappa(\:\!\bm{\cdot}\,;t)$ in~\eqref{eq:kappa} is strictly positive and can be interpreted as a change of numéraire relative to the (dis)utility and corresponding to the coefficients $f_y(\:\!\bm{\cdot}\,;t)$ and $f_z(\:\!\bm{\cdot}\,;t)$.
\end{remark}

\begin{definition}[Adjoint equation/process of first order]
\label{def:p}
We call the \emph{adjoint equation of first order associated with $\rounds{\bar{\u}(\:\!\bm{\cdot}\:\!),\bar{X}(\:\!\bm{\cdot}\:\!),\bar{Y}(\:\!\bm{\cdot}\,;t),\bar{Z}(\:\!\bm{\cdot}\,;t)}$} the linear BSDE
\begin{equation}
\label{eq:p}
\begin{sistema}
d\:\!p(s\:\!;t) = -\:\!g(s,p(s\:\!;t),q(s\:\!;t);t)\:\!ds + q(s\:\!;t)\:\!dW(s), \\[0.5ex]
p(T\:\!;t) = h_x(\bar{X}(T);t),
\end{sistema}
\end{equation} 
(where $s\in\mathopen{[}t,T\mathclose{]}$ and) where the map
\[
g \colon \mathopen{[}0,T\mathclose{]}{_s} \times \R_p \times \R_q \times \mathopen{[}0,T\mathclose{[}{_t} \to \R
\]
is given, for $s\in\mathopen{[}t,T\mathclose{]}$ (and $\probP$-\:\!a.s.), by
\begin{equation}
\label{eq:g}
g(s,p,q\:\!;t) \eqnot \bracks[\big]{b_x(s) + f_z(s\:\!;t)\sigma_x(s) + f_y(s\:\!;t)}p + \bracks[\big]{\sigma_x(s) + f_z(s\:\!;t)}q + f_x(s\:\!;t)
\end{equation} 
(the dependence on $\omega\in\Omega$ is implicit). We call \emph{adjoint process of first order associated with $\rounds{\bar{\u}(\:\!\bm{\cdot}\:\!),\bar{X}(\:\!\bm{\cdot}\:\!),\bar{Y}(\:\!\bm{\cdot}\,;t),\bar{Z}(\:\!\bm{\cdot}\,;t)}$} any process
\[p(\:\!\bm{\cdot}\,;t)
\]
such that $\rounds{\:\!p(\:\!\bm{\cdot}\,;t),q(\:\!\bm{\cdot}\,;t)}$ is a pair solution of~\eqref{eq:p} that belongs to
\[
\LLped{2}{\Omega;\CC\rounds{\mathopen{[}t,T\mathclose{]};\R}}{\F} \times \LLped{2}{t,T;\R}{\F}.
\]
\end{definition} 

\begin{remark}
The process $p(\:\!\bm{\cdot}\,;t)$, or the independent variable $p$ of $g$ in~\eqref{eq:g} are not to be confused with the summability exponent $p$ (see Definition~\ref{def:admissible}), especially because we might want the latter to take the fixed value $p = 2$.
\end{remark}

\begin{definition}[Adjoint equation/process of second order]
\label{def:P}
We call the \emph{adjoint equation of second order associated with $\rounds{\bar{\u}(\:\!\bm{\cdot}\:\!),\bar{X}(\:\!\bm{\cdot}\:\!),\bar{Y}(\:\!\bm{\cdot}\,;t),\bar{Z}(\:\!\bm{\cdot}\,;t)}$}  the linear BSDE
\begin{equation}
\label{eq:P}
\begin{sistema}
dP(s\:\!;t) = -\:\!G(s,P(s\:\!;t),Q(s\:\!;t);t)\:\!ds + Q(s\:\!;t)\:\!dW(s), \\[0.5ex]
P(T\:\!;t) = h_{xx}(\bar{X}(T);t),
\end{sistema}
\end{equation} 
(where $s\in\mathopen{[}t,T\mathclose{]}$ and) where the map
\[
G \colon \mathopen{[}0,T\mathclose{]}{_s} \times \R_P \times \R_Q \times \mathopen{[}0,T\mathclose{[}{_t} \to \R
\]
is given, for $s\in\mathopen{[}t,T\mathclose{]}$ (and $\probP$-\:\!a.s.), by
\begin{multline}
\label{eq:G}
 G(s,P,Q\:\!;t) \eqnot \bracks*{2b_x(s) + {\sigma_x(s)}^2\:\!\! + 2f_z(s\:\!;t)\sigma_x(s) + f_y(s\:\!;t)}\:\!\!P \\[0.75ex]
	+ \bracks[\big]{2\sigma_x(s) + f_z(s\:\!;t)}Q+\:\! b_{xx}(s)p(s\:\!;t) + \sigma_{xx}(s)\bracks[\big]{f_z(s\:\!;t)p(s\:\!;t) + q(s\:\!;t)} \\[0.75ex]
		+ \big{(}1,p(s\:\!;t),\sigma_x(s)p(s\:\!;t) + q(s\:\!;t)\big{)} \cdot \D^2 f(s\:\!;t) \cdot \:\!\! {\big{(}1,p(s\:\!;t),\sigma_x(s)p(s\:\!;t) + q(s\:\!;t)\big{)}}^{\:\!\!\T}\:\!\!
\end{multline} 
(the dependence on $\omega\in\Omega$ is implicit). We call \emph{adjoint process of second order associated with $\rounds{\bar{\u}(\:\!\bm{\cdot}\:\!),\bar{X}(\:\!\bm{\cdot}\:\!),\bar{Y}(\:\!\bm{\cdot}\,;t),\bar{Z}(\:\!\bm{\cdot}\,;t)}$} any process
\[
P(\:\!\bm{\cdot}\,;t)
\]
such that $\rounds{P(\:\!\bm{\cdot}\,;t),Q(\:\!\bm{\cdot}\,;t)}$ is a pair solution of~\eqref{eq:P} that belongs to
\[
\LLped{2}{\Omega;\CC\rounds{\mathopen{[}t,T\mathclose{]};\R}}{\F} \times \LLped{2}{t,T;\R}{\F}.
\]
\end{definition} 

For the adjoint equations~\eqref{eq:p}  and~\eqref{eq:P} of Definitions~\ref{def:p} and~\ref{def:P}, respectively, the following result regarding existence, uniqueness, and regularity holds (see~\cite{hu17}).

\begin{proposition}
\label{prop:pP}
There exist unique pair solutions $\rounds{\:\!p(\:\!\bm{\cdot}\,;t),q(\:\!\bm{\cdot}\,;t)}$ of~\eqref{eq:p} and $\rounds{P(\:\!\bm{\cdot}\,;t),Q(\:\!\bm{\cdot}\,;t)}$ of~\eqref{eq:P} such that, for any $k\in\mathopen{[}1,\infty\mathclose{[}$,
\[
\E\:\!\!\sup_{s\in\mathopen{[}t,T\mathclose{]}} \mathopen{\Big{[}}{\abs*{\:\!p(s\:\!;t)}}^{2k}\:\!\! + {\abs*{P(s\:\!;t)}}^{2k}\mathclose{\Big{]}} + \:\! \E\:\!\!{\mathopen{\bigg{(}}\int_t^T\:\!\! \mathopen{\Big{[}}{\abs*{\:\!q(s\:\!;t)}}^2\:\!\! + {\abs*{Q(s\:\!;t)}}^2\mathclose{\Big{]}}\:\!ds\mathclose{\bigg{)}}}^{\:\!\!\!k}\! < \infty.
\]
\end{proposition} 

\begin{remark}
\label{rem:pPdeterministic}
$p(t\:\!;t)$ and $P(t\:\!;t)$ are deterministic constants (for similar reasons to those given in Remark~\ref{rem:deterministic}), while it is not possible to say the same in general about $q(t\:\!;t)$ and $Q(t\:\!;t)$.
\end{remark} 

\begin{definition}[Generalized Hamiltonian function of first order]
\label{def:H}
We call the \emph{generalized Hamiltonian (function) of first order associated with $b,\sigma,f$} the map
\[
H \colon \mathopen{[}0,T\mathclose{]}{_s} \times \R_x \times U_{\u} \times \R_y \times \R_z \times \R_p \times \R_q \times \mathopen{[}0,T\mathclose{[}{_t} \times I_{\bar{x}} \times U_{\bar{\u}} \to \R
\]
given by
\begin{multline}
\label{eq:H}
H(s,x,\u,y,z,p,q\:\!;t,\bar{x},\bar{\u}) \doteq p\:\!b(s,x,\u) + q\:\!\sigma(s,x,\u) \\[0.5ex]
	 + f(s,x,\u,y,z + p\:\!\bracks{\sigma(s,x,\u) - \sigma(s,\bar{x},\bar{\u})};t).
\end{multline} 
\end{definition} 

\begin{definition}[Generalized Hamiltonian function of second order]
\label{def:HH}
We call the \emph{generalized Hamiltonian (function) of second order associated with $b,\sigma,f$} the map
\[
\HH \colon \mathopen{[}0,T\mathclose{]}{_s} \times \R_x \times U_{\u} \times \R_y \times \R_z \times \R_p \times \R_q \times \R_P \times \mathopen{[}0,T\mathclose{[}{_t} \times I_{\bar{x}} \times U_{\bar{\u}} \to \R
\]
given by
\begin{equation}
\label{eq:HH}
\HH(s,x,\u,y,z,p,q,P\:\!;t,\bar{x},\bar{\u}) \doteq H(s,x,\u,y,z,p,q\:\!;t,\bar{x},\bar{\u}) + \textstyle{\frac{1}{2}}\:\!P\;\!{\bracks{\sigma(s,x,\u) - \sigma(s,\bar{x},\bar{\u})}}^2\:\!\!.
\end{equation} 
\end{definition} 

\begin{notation}
\label{not:HH}
With respect to $\rounds{\bar{\u}(\:\!\bm{\cdot}\:\!),\bar{X}(\:\!\bm{\cdot}\:\!),\bar{Y}(\:\!\bm{\cdot}\,;t),\bar{Z}(\:\!\bm{\cdot}\,;t)}$ and the corresponding adjoint processes $p(\:\!\bm{\cdot}\,;t)$ and $P(\:\!\bm{\cdot}\,;t)$ as in Definitions~\ref{def:p} and~\ref{def:P} respectively, and for $s\in\mathopen{[}t,T\mathclose{]}$ (and $\probP$-\:\!a.s.), we write
\begin{equation*}
 \HH(s\:\!;t) \eqnot \HH(s,\bar{X}(s),\bar{\u}(s),\bar{Y}(s\:\!;t),\bar{Z}(s\:\!;t),p(s\:\!;t),q(s\:\!;t),P(s\:\!;t);t,\bar{X}(s),\bar{\u}(s))
\end{equation*}
and, for $\u \in U$,
\begin{equation*}
\delta\HH(s\:\!;t,\u) \eqnot \HH(s,\bar{X}(s),\u,\bar{Y}(s\:\!;t),\bar{Z}(s\:\!;t),p(s\:\!;t),q(s\:\!;t),P(s\:\!;t);t,\bar{X}(s),\bar{\u}(s)) - \HH(s\:\!;t).
\end{equation*}
\end{notation} 

Regarding the following, let us keep in mind Notation~\ref{not:HH}.

\begin{remark}
For $s\in\mathopen{[}t,T\mathclose{]}$ and $\u \in U$, $\HH(s\:\!;t)$ and $\delta\HH(s\:\!;t,\u)$ belong to $\Lped{1}{\Omega;\R}{s}$ and, moreover,
\begin{equation}
\label{eq:deltaHH0}
\delta\HH(s\:\!;t,\bar{\u}(s)) = 0.
\end{equation} 
\end{remark}

The key result is the following lemma, which we will prove in~\ref{subsec:lemmaliminf}.

\begin{lemma}
\label{lemma:liminf}
Fix $\u(\:\!\bm{\cdot}\:\!)\in\UUU\bracks{t,T}$. For $\epsilon\in\mathopen{]}0,T-t\mathclose{[}$ and $s\in\mathopen{[}t,T\mathclose{]}$, let $E_t^{\:\!\epsilon} \equiv E_{t,s}^{\:\!\epsilon}$\:\!\! be given by
\begin{equation}
\label{eq:Etseps}
E_t^{\:\!\epsilon} \eqnot
\begin{cases}
\mathopen{[}s,s+\epsilon\mathclose{]}, & \text{if $s<T$,} \\[0.5ex]
\mathopen{[}T-\epsilon,T\mathclose{]}, & \text{if $s=T$,}
\end{cases}
\end{equation} 
and let $\bar{\u}^\epsilon(\:\!\bm{\cdot}\:\!)$ be the spike variation of $\bar{\u}(\:\!\bm{\cdot}\:\!)$ w.r.t. $\u(\:\!\bm{\cdot}\:\!)$ and $E_t^{\:\!\epsilon}$\:\!\!. Then, for any $s\in\mathopen{[}t,T\mathclose{]}$, the $\liminf_{\epsilon\downarrow0}$ as in~\eqref{eq:liminf} is an actual limit and  takes the form
\begin{equation}
\label{eq:lemmaliminf}
\liminf_{\epsilon\downarrow0}\:\!\frac{J(\bar{\u}^\epsilon(\:\!\bm{\cdot}\:\!);t,x) - J(\bar{\u}(\:\!\bm{\cdot}\:\!);t,x)}{\epsilon} = \E\bracks[\big]{\kappa(s\:\!;t)\:\!\delta\HH(s\:\!;t,\u(s))}.
\end{equation} 
\end{lemma} 

\begin{remark}
As we will understand shortly, instead of~\eqref{eq:Etseps}, we could take, among other possibilities,
\[
E_t^{\:\!\epsilon} \eqnot
\begin{cases}
\mathopen{[}s,s+\epsilon\mathclose{[}, & \text{if $s<T$,} \\[0.5ex]
\mathopen{]}T-\epsilon,T\mathclose{]}, & \text{if $s=T$}.
\end{cases}
\]
\end{remark}

\begin{corollary}[Sufficient conditions]
\label{cor:liminf}
Suppose there exists a measurable map $\Py \colon \mathopen{[}0,T\mathclose{]} \times I \to U$ such that, for any $s\in\mathopen{[}t,T\mathclose{]}$ (and $\probP$-\:\!a.s.),
\begin{equation}
\label{eq:uPy}
\bar{\u}(s) = \Py(s,\bar{X}(s))
\end{equation} 
and suppose that, for any $\u \in U$ (and $\probP$-\:\!a.s.),
\begin{equation}
\label{eq:deltaHH}
\delta\HH(t\:\!;t,\u) \geq 0.
\end{equation} 
Then $\Py$ is an equilibrium policy, i.e.,
\[
\rounds{\bar{\u}(\:\!\bm{\cdot}\:\!),\bar{X}(\:\!\bm{\cdot}\:\!),\bar{Y}(\:\!\bm{\cdot}\,;t),\bar{Z}(\:\!\bm{\cdot}\,;t)}
\]
is an equilibrium 4-tuple.
\end{corollary} 
\begin{proof}
By~\eqref{eq:lemmaliminf} of Lemma~\ref{lemma:liminf}, with $E_t^{\:\!\epsilon}$\:\!\! as in~\eqref{eq:Eteps}, that is,
\[
s = t,
\]
the inequality~\eqref{eq:liminf} holds (see also Definition~\ref{def:kappa}). 
\end{proof}

\begin{remark}
\label{rem:corliminf}
Regarding Corollary~\ref{cor:liminf}, we point out that, \emph{if} $\bar{Z}(t\:\!;t)$ and $q(t\:\!;t)$ are assumed to be deterministic constants, then the condition~\eqref{eq:deltaHH} is equivalent to
\begin{equation}
\label{eq:Pyargmin}
\Py(t,x) \in \argmin_{\u \in U} \HH(t,x,\u,\bar{Y}(t\:\!;t),\bar{Z}(t\:\!;t),p(t\:\!;t),q(t\:\!;t),P(t\:\!;t);t,x,\Py(t,x))
\end{equation} 
because, by~\eqref{eq:uPy},
\[
\Py(t,x) = \bar{\u}(t).
\]
See Section~\ref{sec:applicationuncon} to figure out that assuming $\bar{Z}(t\:\!;t)$ and $q(t\:\!;t)$ as deterministic may be superfluous or circumventable in concrete applications.
\end{remark} 

We are finally ready to present the first of our main results (see also Corollary~\ref{cor:liminf} and Remark~\ref{rem:corliminf}).

\begin{theorem}[Maximum principle] 
\label{th:maxprinc}
Suppose there exists a measurable map $\Py \colon \mathopen{[}0,T\mathclose{]} \times I \to U$ such that, for any $s\in\mathopen{[}t,T\mathclose{]}$ (and $\probP$-\:\!a.s.),
\[
\bar{\u}(s) = \Py(s,\bar{X}(s)).
\]
Then the following three conditions are equivalent.
\begin{itemize}[leftmargin=*]
	\item[\textbf{1.}] $\Py$ is an equilibrium policy, i.e., $\rounds{\bar{\u}(\:\!\bm{\cdot}\:\!),\bar{X}(\:\!\bm{\cdot}\:\!),\bar{Y}(\:\!\bm{\cdot}\,;t),\bar{Z}(\:\!\bm{\cdot}\,;t)}$ is an equilibrium 4-tuple.
	\item[\textbf{2.}] For any $\u \in U$ (and $\probP$-\:\!a.s.),
	\[
	\delta\HH(t\:\!;t,\u) \geq 0.
	\]
	\item[\textbf{3.}] If $\bar{Z}(t\:\!;t)$ and $q(t\:\!;t)$ are deterministic constants, then
		\[
		\Py(t,x) \in \argmin_{\:\!\u \in U} \HH(t,x,\u,\bar{Y}(t\:\!;t),\bar{Z}(t\:\!;t),p(t\:\!;t),q(t\:\!;t),P(t\:\!;t);t,x,\Py(t,x)).
		\]
\end{itemize}
\end{theorem} 
\begin{proof}
In light of what we saw in Corollary~\ref{cor:liminf} and Remark~\ref{rem:corliminf}, we just need to show that
\[
\textbf{1 $\Rightarrow$ 2}
\]
(necessary conditions, we would now say). To this end suppose, by contradiction, that there exist $t^\ast\:\!\!\in\mathopen{[}0,T\mathclose{[}$, $\u^{\:\!\!\ast}\:\!\! \in U$, and $\NN\in\FF$ with $\probP\bracks{\:\!\NN\:\!} > 0$, such that
\[
\delta\HH(t^\ast;t^\ast\:\!\!,\u^{\:\!\!\ast}) < 0
\]
on $\NN$. Then any $\u(\:\!\bm{\cdot}\:\!)\in\UUU\bracks{t^\ast\:\!\!,T}$ such that, $\probP$-\:\!a.s.,
\[
\u(t^\ast) =
\begin{cases}
\bar{\u}(t^\ast), & \text{on $\Omega\setminus\NN$,} \\[0.25ex]
\u^{\:\!\!\ast\:\!\!}, & \text{on $\NN$,}
\end{cases}
\]
(e.g., the trivial one) satisfies
\[
\delta\HH(t^\ast;t^\ast\:\!\!,\u(t^\ast)) = \delta\HH(t^\ast;t^\ast\:\!\!,\u^{\:\!\!\ast})\:\!\1_\NN
\]
(see also~\eqref{eq:deltaHH0}), and therefore, since $\probP\bracks{\:\!\NN\:\!} > 0$,
\[
\E\bracks[\big]{\delta\HH(t^\ast;t^\ast\:\!\!,\u(t^\ast))} < 0,
\]
which is a contradiction (see also~\eqref{eq:lemmaliminf}). 
\end{proof}

Under appropriate assumptions on our coefficients, we can replace $\HH$ with $H$ in Theorem~\ref{th:maxprinc} thus obtaining the following result, which will be used in Section~\ref{sec:applicationuncon}. Let us keep in mind also Definition~\ref{def:P}.

\begin{corollary}
\label{cor:maxprinc}
Suppose there exists a measurable map $\Py \colon \mathopen{[}0,T\mathclose{]} \times I \to U$ such that, for any $s\in\mathopen{[}t,T\mathclose{]}$ (and $\probP$-\:\!a.s.),
\[
\bar{\u}(s) = \Py(s,\bar{X}(s))
\]
and suppose that $\bar{Z}(t\:\!;t)$ and $q(t\:\!;t)$ are deterministic constants. If $h(\:\!\bm{\cdot}\,;t)$ is convex and
\[
G(\:\!\bm{\cdot}\,,0,0\:\!;t) \geq 0,
\]
then the following two conditions are equivalent.
\begin{itemize}[leftmargin=*]
	\item[\textbf{1.}] $\Py$ is an equilibrium policy, i.e., $\rounds{\bar{\u}(\:\!\bm{\cdot}\:\!),\bar{X}(\:\!\bm{\cdot}\:\!),\bar{Y}(\:\!\bm{\cdot}\,;t),\bar{Z}(\:\!\bm{\cdot}\,;t)}$ is an equilibrium 4-tuple.
	\item[\textbf{2.}] $\Py(t,x) \in \argmin_{\:\!\u \in U} H(t,x,\u,\bar{Y}(t\:\!;t),\bar{Z}(t\:\!;t),p(t\:\!;t),q(t\:\!;t);t,x,\Py(t,x))$.
\end{itemize}
\end{corollary} 
\begin{proof}
This is a direct consequence of Theorem~\ref{th:maxprinc}: indeed, comparing (\eqref{eq:P} and)~\eqref{eq:G} of Definition~\ref{def:P} with~\eqref{eq:linearBSDEs1} of Proposition~\ref{prop:linearBSDEs} (Section~\ref{sec:preliminary}), we deduce that
\[
P(\:\!\bm{\cdot}\,;t) \geq 0
\]
by~\eqref{eq:comparison} of Remark~\ref{rem:comparison} (where $\xi_{\:\!t} \equiv h_{xx}(\bar{X}(T);t)$, $\alpha(\:\!\bm{\cdot}\,;t) \equiv G(\:\!\bm{\cdot}\,,0,0\:\!;t)$ and $\Xi(\:\!\bm{\cdot}\,;t) \equiv P(\:\!\bm{\cdot}\,;t)$) and so the term of $\HH$ that depends on it, namely,
\[
\textstyle{\frac{1}{2}}\:\!P(t\:\!;t)\:\!{\bracks{\sigma(t,x,\u) - \sigma(t,\bar{x},\bar{\u}(t))}}^2\:\!\!
\]
is superfluous in calculating the minimum (as in~\eqref{eq:Pyargmin}). 
\end{proof}

It is somewhat standard to define a multidimensional Wiener process, some regular coefficients, an open pluri-interval state domain $I$, a recursive stochastic control problem, the adjoint equations/processes of first and second order (associated with any admissible 4\:\!-tuple) and the generalized Hamiltonian of second order (associated with the coefficients) in such a way that, if the (dis)utility functional $J(\:\!\bm{\cdot}\,;t,x)\colon\UUU\bracks{t,T}\to\R$ changes from~\eqref{eq:utilityfun} to
\[
J(\u(\:\!\bm{\cdot}\:\!);t,x) \doteq \gamma(Y(t\:\!;t);t)
\]
($t\in\mathopen{[}0,T\mathclose{[}$ and $x \in I$), where
\[
\gamma \colon \R_y^m\:\!\! \times \mathopen{[}0,T\mathclose{[}{_t} \to \R
\]
is of class $C^1$\:\!\! w.r.t. the variable $y\in\R^m$\:\!\! ($m\in\N^*$\:\!\!), then (Lemma~\ref{lemma:liminf} and) Theorem~\ref{th:maxprinc} can be extended.

\begin{theorem} 
\label{th:maxprincmulti}
Suppose there exists a measurable map $\Py \colon \mathopen{[}0,T\mathclose{]} \times I \to U$ such that, for any $s\in\mathopen{[}t,T\mathclose{]}$ (and $\probP$-\:\!a.s.), $\bar{\u}(s) = \Py(s,\bar{X}(s))$. Then the following three conditions are equivalent.
\begin{itemize}[leftmargin=*]
	\item[\textbf{1.}] $\Py$ is an equilibrium policy, i.e., $\rounds{\bar{\u}(\:\!\bm{\cdot}\:\!),\bar{X}(\:\!\bm{\cdot}\:\!),\bar{Y}(\:\!\bm{\cdot}\,;t),\bar{Z}(\:\!\bm{\cdot}\,;t)}$ is an equilibrium 4-tuple.
	\item[\textbf{2.}] For any $\u \in U$ (and $\probP$-\:\!a.s.),
	\[
	{\big{(}}\D_y\;\!\gamma(\bar{Y}(t\:\!;t);t){\big{)}}^{\!\T\:\!\!} \:\!\! \cdot \delta\HH(t\:\!;t,\u) \geq 0.
	\]
	\item[\textbf{3.}] If $\bar{Z}(t\:\!;t)$ and $q(t\:\!;t)$ are deterministic constants, then
		\[
		\Py(t,x) \in \argmin_{\:\!\u \in U} {\big{(}}\D_y\;\!\gamma(\bar{Y}(t\:\!;t);t){\big{)}}^{\!\T\:\!\!} \:\!\! \cdot \HH(t,x,\u,\bar{Y}(t\:\!;t),\bar{Z}(t\:\!;t),p(t\:\!;t),q(t\:\!;t),P(t\:\!;t);t,x,\Py(t,x)).
		\]
\end{itemize}
\end{theorem} 

It may be useful to consult~\cite{hu17} on this topic, especially regarding the dimensional changes related to the adjoint processes and the generalized Hamiltonian.


\section{Application to portfolio management} 
\label{sec:applicationuncon}

In this section, we investigate the Merton portfolio management problem with recursive utility in the special case of non-exponential discounting and CRRA preferences. In the search for a concrete equilibrium policy $\Py$, we first fit the sufficient conditions of Theorem~\ref{th:maxprinc} to the portfolio problem. Then, we set up a sort of generalized HJB equation for the equilibrium policy. Under precise assumptions, we characterize the equilibrium policy by imposing a suitable \emph{ansatz} for the value function associated with the problem itself. It turns out that the equilibrium policy is a solution of a partial differential equation whose existence and uniqueness are also investigated. In this regard see, e.g.,~\cite{ekelandpirvu08} and~\cite{ekelandmbodjipirvu12}. 




In a Black\:\!–\:\!Scholes financial market model, let $S_0 = (S_0(s))_{s\in\mathopen{[}0,T\mathclose{]}}$ be the value of a (deterministic) saving account, or \emph{bond}, which accrues interest at a given constant risk-free rate
\[
r \in \mathopen{]}0,\infty\mathclose{[}
\]
(with $S_0(0) \in \mathopen{]}0,\infty\mathclose{[}$ exogenously specified): that is, for $s\in\mathopen{[}0,T\mathclose{]}$,
\[
dS_0(s) = r\:\!S_0(s)\:\!ds,
\]
i.e., $S_0(s) = S_0(0)\:\!e^{rs}$\:\!\!. Let $S = (S(s))_{s\in\mathopen{[}0,T\mathclose{]}}$ be the price of a risky asset, or \emph{stock}, which evolves as a geometric, or exponential, Brownian motion $W$\:\!\!-\:\!based (with $S(0) \in \mathopen{]}0,\infty\mathclose{[}$ exogenously specified): that is, for $s\in\mathopen{[}0,T\mathclose{]}$,
\[
dS(s) = S(s)\bracks[\big]{\rho\;\!ds + \sigma\:\!dW(s)},
\]
i.e., explicitly,
\[
S(s) = S(0)\exp\braces[\Big]{\:\!\!(\rho-\sigma^2\:\!\!/2)\:\!s + \sigma\:\!W(s)\:\!\!},
\]
where
\[
\rho\in\mathopen{]}r,\infty\mathclose{[}
\]
is the mean rate of return by the stock (or appreciation rate of the stock) and
\[
\sigma\in\mathopen{]}0,\infty\mathclose{[}
\]
is the volatility of the price (both of which are constant scalars). The generalization to the case of multiple stocks would basically be a matter of formalization.

We assume that the above market is complete (there are no sources of randomness other than the stock) and we denote the excess return, on average, by investment in the stock as
\[
\mu \eqnot \rho - r\:\!.
\]

A decision-maker within this market, or agent for short, is assumed to invest her/his wealth
\[
X(\:\!\bm{\cdot}\:\!),
\]
value process of the portfolio, in the stock and the bond and to consume continuously over time (on $\mathopen{[}0,T\mathclose{]}$), depending dynamically on the current wealth itself.

Therefore, for $s\in\mathopen{[}0,T\mathclose{]}$ (and $\probP$-\:\!a.s.), let
\[
\zeta(s,X(s)) \in \R
\]
be the proportion of current wealth $X(s) \in \mathopen{]}0,\infty\mathclose{[}$ invested in the stock at time $s$ (with a sign), and
\[
c(s,X(s)) \in \mathopen{[}0,\infty\mathclose{[}
\]
be the proportion of $X(s)$ consumed at time $s$. Then $n = 2$, the control domain becomes
\begin{equation}
\label{eq:applicationdomain}
U = \R{_\zeta} \times \mathopen{[}0,\infty\mathclose{[}{_c},
\end{equation} 
and any admissible control $\u(\:\!\bm{\cdot}\:\!)\in\UUU\bracks{t,T}$ is a portfolio strategy that can be written as
\begin{equation}
\label{eq:feedback}
\u(\:\!\bm{\cdot}\:\!) = \rounds{\zeta(\:\!\bm{\cdot}\,,X(\:\!\bm{\cdot}\:\!)),c(\:\!\bm{\cdot}\,,X(\:\!\bm{\cdot}\:\!))},
\end{equation} 
and which we can refer to as \emph{investment-consumption policy} or \emph{trading-consumption policy} (see also Definition~\ref{def:admissible}). We remark that any $\u(\:\!\bm{\cdot}\:\!)$ in~\eqref{eq:feedback} depends on $\omega\in\Omega$ only through $X(\:\!\bm{\cdot}\:\!)$ and that appears in feedback or closed-loop form w.r.t. $X(\:\!\bm{\cdot}\:\!)$ itself.

\begin{remark}
$1 - \zeta(\:\!\bm{\cdot}\,,X(\:\!\bm{\cdot}\:\!))$ coincides with the proportion of $X(\:\!\bm{\cdot}\:\!)$ invested in the bond.
\end{remark}

The equation describing the dynamics of the wealth process $X(\:\!\bm{\cdot}\:\!) \equiv X^{\u}(\:\!\bm{\cdot}\:\!)$ is the controlled (forward) stochastic differential equation on $\mathopen{[}0,T\mathclose{]}\times\Omega$, in It\^o differential form, given by
\begin{equation}
\label{eq:wealth}
dX(s) = X(s)\bracks*{\mathopen{\big{(}}r + \mu\:\!\zeta(s,X(s)) - c(s,X(s))\mathclose{\big{)}}\:\!ds + \sigma\:\!\zeta(s,X(s))\:\!dW(s)},
\end{equation} 
which, for any possible initial wealth $X(0) \in \mathopen{]}0,\infty\mathclose{[}$, admits an unique positive solution
\[
X(\:\!\bm{\cdot}\:\!) \in \LLped{2}{\Omega;\CC\rounds{\mathopen{[}0,T\mathclose{]};\R}}{\F}.
\]
In particular, the usual self\:\!-\:\!financing condition is satisfied: namely, the variation in wealth over time is due exclusively to profits and losses from investing in the stock and from consumption (there is no cashflow coming in or out).

Now, we assume that the agent derives utility from intertemporal consumption
\[
c(\:\!\bm{\cdot}\,,X(\:\!\bm{\cdot}\:\!))X(\:\!\bm{\cdot}\:\!)
\]
and final wealth $X(T)$, which she/he tries to optimize by minimizing, at least in a weak sense, a discounted expectation---or, rather, a recursive (dis)utility---involving (dis)utility functions.

Therefore, let $\upsilon(\:\!\bm{\cdot}\:\!)$ and $\hat{\upsilon}(\:\!\bm{\cdot}\:\!)$ be two scalar functions of a real variable that satisfy the classical Uzawa--Inada conditions (utility functions, in fact): i.e.,
\[
\upsilon\colon\mathopen{[}0,\infty\mathclose{[}\to\mathopen{[}0,\infty\mathclose{[}
\]
is of class $C^2$\:\!\! and strictly increasing such that, for any $x\in\mathopen{]}0,\infty\mathclose{[}$, $\upsilon''(x) < 0$ (thus, $\upsilon(\:\!\bm{\cdot}\:\!)$ is strictly concave on $\mathopen{]}0,\infty\mathclose{[}$), with $\upsilon(0) = 0$ and
\[
\lim_{\:\!x\downarrow0} \upsilon'(x) = \infty, \qquad \lim_{\:\!x\uparrow\infty} \upsilon'(x) = 0
\]
(the same for $\hat{\upsilon}(\:\!\bm{\cdot}\:\!)$).

\begin{remark}
\label{rem:Upsilon}
The marginal disutility function $\upsilon'\colon\mathopen{[}0,\infty\mathclose{[}\to\mathopen{]}0,\infty\mathclose{[}$ is of class $C^1$\:\!\! and bijective outside the origin, and has an inverse function
\[
(\upsilon')^{-1}\colon\mathopen{]}0,\infty\mathclose{[}\to\mathopen{]}0,\infty\mathclose{[}
\]
that is continuous and strictly decreasing with, again,
\[
\lim_{\:\!x\downarrow0} (\upsilon')^{-1}(x) = \infty, \qquad \lim_{\:\!x\uparrow\infty} (\upsilon')^{-1}(x) = 0
\]
and, in particular, $(-\:\!\upsilon')^{-1}(\:\!\bm{\cdot}\:\!)$ is a positive function with domain $\mathopen{]}-\infty,0\mathclose{[}$ (the same for $\hat{\upsilon}'(\:\!\bm{\cdot}\:\!)$).
\end{remark} 


\begin{notation}
\label{not:Upsilon}
We write
\[
\Upsilon \eqnot (-\:\!\upsilon')^{-1}\:\!\!.
\]
\end{notation} 

\begin{example}[CRRA preferences]
\label{ex:utilityfunctions}
Regarding $\upsilon(\:\!\bm{\cdot}\:\!)$, we could take $\lambda \in \mathopen{]}0,1\mathclose{[}$ and, for $x\in\mathopen{[}0,\infty\mathclose{[}$,
\[
\upsilon(x) \equiv \upsilon_{\lambda}(x) \doteq x^{\lambda}/\lambda
\]
thus obtaining a constant relative risk aversion equal to
\[
1-\lambda \equiv - \:\! x \:\! \upsilon''(x)/\upsilon'(x)
\]
and furthermore, for any $y \in \mathopen{]}0,\infty\mathclose{[}$,
\[
\Upsilon(- \:\! y) = y^{- \:\! \frac{1}{1 - \lambda}} \:\!\!
\]
(similarly for $\hat{\upsilon}(\:\!\bm{\cdot}\:\!)$).
\end{example} 

Also, for any fixed $t\in\mathopen{[}0,T\mathclose{[}$, let $\hslash(\:\!\bm{\cdot}\,;t)$ and $\hat{\hslash}(\:\!\bm{\cdot}\,;t)$ be two regular discount functions on $\mathopen{[}t,T\mathclose{]}$: i.e.,
\[
\hslash(\:\!\bm{\cdot}\,;t)\colon\mathopen{[}t,T\mathclose{]}\to\mathopen{]}0,\infty\mathclose{[}
\]
of class $C^1$\:\!\! with
\[
\hslash(t\:\!;t) = 1
\]
(the same for $\hat{\hslash}(\:\!\bm{\cdot}\,;t)$). 	

\begin{remark}
The corresponding (psychological) discount rate
\[
- \:\! \hslash'(\:\!\bm{\cdot}\,;t)/\hslash(\:\!\bm{\cdot}\,;t),
\]
i.e., the rate of return used to discount future cashflows back to their present value, should be regarded as a monotonic function (the same for $\hat{\hslash}(\:\!\bm{\cdot}\,;t)$).
\end{remark}

We must emphasize here that time\:\!-inconsistency arises from the way the future is discounted, and not from a change in preferences (which are the same for all the decision-makers).

\begin{example}[Non\:\!-\:\!exponential discounting]
\label{ex:discountfunctions}
Regarding $\hslash(\:\!\bm{\cdot}\,;t)$, we could take a common regular discount function $h(\:\!\bm{\cdot}\:\!)$ on $\mathopen{[}0,T-t\mathclose{]}$ and define, for $s\in\mathopen{[}t,T\mathclose{]}$,
\[
\hslash(s\:\!;t) \doteq h(s - t)
\]
where for instance, associated with $K\in\mathopen{]}0,\infty\mathclose{[}$, we could imagine that, for $\tau\in\mathopen{[}0,T-t\mathclose{]}$,
\[
h(\tau) \equiv h_K (\tau) \doteq \frac{1}{1+K\tau}
\]
then have to deal with hyperbolic discounting and non-constant discount rate (similarly for $\hat{\hslash}(\:\!\bm{\cdot}\,;t)$).
\end{example} 

\begin{remark}
A non\:\!-\:\!exponential discount mechanism, such as the one in Example~\ref{ex:discountfunctions}, inevitably leads to time\:\!-inconsistency, and we could also interpret this phenomenon by means of non\:\!-\:\!decomposable financial laws: i.e., given any cashflow $C(\:\!\bm{\cdot}\:\!)$ on $\mathopen{]}0,T\mathclose{]}$, whose utility through $\upsilon(\:\!\bm{\cdot}\:\!)$ must be actualized (from future to present), it does not hold that, for any triplet $t,s,\tau$ with $0 < t < s < \tau \leq T$,
\[
h(\tau-s) \:\! h(s-t) = h(\tau-t)
\]
since, indeed, this is true if and only if $h(s-t) = e^{\:\!s-t}$\:\!\! (as is known and easy to check). In this regard, the diagram below might help.
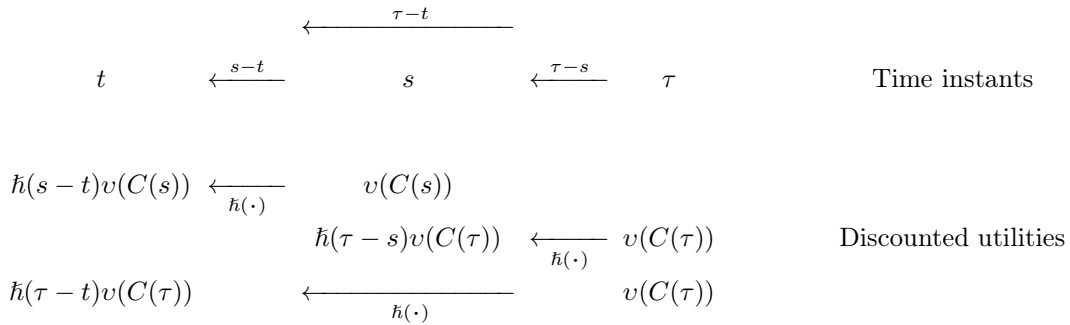
\[
\begin{CD}
			@.											@.	@. \\
			@.	\xleftarrow{\quad\qquad \tau-t \qquad\quad}		@.	@. \\
\text{$t$}	@< s-t <<	\text{$s$}	@< \tau-s <<	\text{$\tau$}		@.	\text{\hspace{10ex} Time instants} \\
			@.											@.	@. \\
\hslash(s-t)\upsilon(C(s))	@<< \hslash(\:\!\bm{\cdot}\:\!) <	\upsilon(C(s))		@. \\
			@.	\hslash(\tau-s)\upsilon(C(\tau))	@<<\hslash(\:\!\bm{\cdot}\:\!)<	\upsilon(C(\tau))		@.	\text{\hspace{10ex} Discounted utilities} \\
\hslash(\tau-t)\upsilon(C(\tau))	@.	\xleftarrow[\quad\qquad \hslash(\:\!\bm{\cdot}\:\!) \qquad\quad]{}		@.	\upsilon(C(\tau))		@. \\
			@.											@.	@.
\end{CD}
\]
\end{remark}

Finally, choosing	
\[
I = \mathopen{]}0,\infty\mathclose{[}
\]
(as the state domain), a fixed $t\in\mathopen{[}0,T\mathclose{[}$, two positive processes
\[
\beta(\:\!\bm{\cdot}\,;t),\gamma(\:\!\bm{\cdot}\,;t) \in \LLped{\infty}{t,T;\R}{\F},
\]
$x \in I$ and $\u(\:\!\bm{\cdot}\:\!) = \rounds{\zeta(\:\!\bm{\cdot}\,,X(\:\!\bm{\cdot}\:\!)),c(\:\!\bm{\cdot}\,,X(\:\!\bm{\cdot}\:\!))} \in \UUU\bracks{t,T}$, we consider the recursive stochastic control problem having the state equation as in~\eqref{eq:wealth} and the recursive (dis)utility system according to the classic Uzawa type (see also Definition~\ref{def:fbsde}): that is,
\begin{equation}
\label{eq:applicationuncon}
\begin{sistema}
dX(s) = X(s)\bracks*{\mathopen{\big{(}}r + \mu\:\!\zeta(s,X(s)) - c(s,X(s))\mathclose{\big{)}}\:\!ds + \sigma\:\!\zeta(s,X(s))\:\!dW(s)}, \\[1ex]
dY(s\:\!;t) = -\:\!\hslash(s\:\!;t)\bracks[\big]{-\:\!\upsilon(c(s,X(s))X(s)) - \beta(s\:\!;t)Y(s\:\!;t) - \gamma(s\:\!;t)Z(s\:\!;t)}\:\!ds \\[1ex]
	\hphantom{dY(s\:\!;t) = \:\!} + Z(s\:\!;t)\:\!dW(s), \\[1ex]
X(t) = x,\quad Y(T\:\!;t) = -\:\!\hat{\hslash}(T;t)\hat{\upsilon}(X(T)),
\end{sistema}
\end{equation} 
(where $s\in\mathopen{[}t,T\mathclose{]}$). See, e.g.,~\cite{elkarouipengquenez97}.

\begin{remark}
By virtue of Remark~\ref{rem:comparison}, also fundamental below, $Y(\:\!\bm{\cdot}\,;t) < 0$ ($\probP$-\:\!a.s.).
\end{remark}

It is quite simple to verify that the sufficient conditions of Theorem~\ref{th:maxprinc}, or rather of Corollary~\ref{cor:maxprinc}, can be expressed as follows and that, indeed, are pratically based on the calculation of the gradient w.r.t. the variable $\u = \rounds{\zeta,c} \in \R \times \mathopen{]}0,\infty\mathclose{[}$ of the generalized Hamiltonian of first order $H$, associated with the coefficients maps of~\eqref{eq:applicationuncon}, where $s = t$, $\bar{x} = x$ and $\bar{\u} = \rounds{\bar{\zeta},\bar{c}}$ (see Definition~\ref{def:H}): i.e.,
\begin{multline}
\label{eq:applicationH}
H(t,x,\rounds{\zeta,c},y,z,p,q\:\!;t,x,\rounds{\bar{\zeta},\bar{c}}) = x\:\!p\:\!\rounds{\mu\:\!\zeta - c} + x\:\!\sigma\:\!q\:\!\zeta	\\[0.25ex]
	- \braces[\Big]{\:\!\!\upsilon(x\:\!c) + y\:\!\beta(t\:\!;t) + \gamma(t\:\!;t)\bracks[\big]{z + x\:\!\sigma\:\!p\:\!\rounds{\zeta - \bar{\zeta}\:\!}}\:\!\!}.
\end{multline} 

\begin{remark}
The gradient of $H$ w.r.t. $\rounds{\zeta,c}$ does not depend on the variable $z$ (see~\eqref{eq:applicationH}) and, therefore, it will not be necessary to assume that $\bar{Z}(t\:\!;t)$ is a deterministic constant in the following theorem.
\end{remark}

\begin{theorem}
\label{th:applicationuncon}
Suppose there exists a measurable map
\[
\Py = \rounds{\Pi_1,\Pi_2} \colon \mathopen{[}0,T\mathclose{]}{_s} \times \mathopen{]}0,\infty\mathclose{[}{_x} \to \R \times \mathopen{[}0,\infty\mathclose{[}
\]
of class $C^2$\:\!\! w.r.t. the variable $x$ with bounded first and second derivatives, even if multiplied by the identity function, such that if, for any $t\in\mathopen{[}0,T\mathclose{[}$ and $x \in \mathopen{]}0,\infty\mathclose{[}$, $\bar{X}(\:\!\bm{\cdot}\:\!)$ is the solution of the FSDE 
\begin{equation}
\label{eq:applicationX}
\begin{sistema}
dX(s) = X(s)\bracks[\Big]{\mathopen{\big{(}}r + \mu\:\!\Pi_1(s,X(s)) - \Pi_2(s,X(s))\mathclose{\big{)}}\:\!ds + \sigma\:\!\Pi_1(s,X(s))\:\!dW(s)}, \\[0.75ex]
X(t) = x,
\end{sistema}
\end{equation} 
(where $s\in\mathopen{[}t,T\mathclose{]}$), then $\Py(\:\!\bm{\cdot}\,,\bar{X}(\:\!\bm{\cdot}\:\!)) \in \UUU\bracks{t,T}$ and there exists a pair solution
\[
\rounds{p(\:\!\bm{\cdot}\,;t),q(\:\!\bm{\cdot}\,;t)} \in \LLped{2}{\Omega;\CC\rounds{\mathopen{[}t,T\mathclose{]};\R}}{\F} \times \LLped{2}{t,T;\R}{\F}
\]
of the BSDE
\begin{equation}
\label{eq:applicationuncon0}
\begin{sistema}
d\:\!p(s\:\!;t) = - \, \mathopen{\bigg{\lbrace}} \:\!\! \mathopen{\Big{\lbrace}} r + \mu\bracks*{\Pi_1(s,\bar{X}(s)) + \bar{X}(s)\frac{\partial}{\partial x}\Pi_1(s,\bar{X}(s))} - \big{\lbrack}\Pi_2(s,\bar{X}(s)) \\[0.75ex]
		\hphantom{d\:\!p(s\:\!;t) = - \, \mathopen{\bigg{\lbrace}} \:\!} + \bar{X}(s)\frac{\partial}{\partial x}\Pi_2(s,\bar{X}(s))\big{\rbrack} - \sigma\:\!\hslash(s\:\!;t)\gamma(s\:\!;t)\big{\lbrack}\Pi_1(s,\bar{X}(s)) \\[1ex]
		\hphantom{d\:\!p(s\:\!;t) = - \, \mathopen{\bigg{\lbrace}} \:\!} + \bar{X}(s)\frac{\partial}{\partial x}\Pi_1(s,\bar{X}(s))\big{\rbrack} - \hslash(s\:\!;t)\beta(s\:\!;t) \mathclose{\Big{\rbrace}}\:\!p(s\:\!;t) \\[1ex]
	\hphantom{d\:\!p(s\:\!;t) = - \, \:\!} + \braces[\Big]{\sigma\bracks*{\Pi_1(s,\bar{X}(s)) + \bar{X}(s)\frac{\partial}{\partial x}\Pi_1(s,\bar{X}(s))} - \hslash(s\:\!;t)\gamma(s\:\!;t)}\:\!q(s\:\!;t) \\[2ex]
	\hphantom{d\:\!p(s\:\!;t) = - \, \:\!} - \hslash(s\:\!;t)\:\!\upsilon'(\Pi_2(s,\bar{X}(s))\bar{X}(s))\big{\lbrack}\Pi_2(s,\bar{X}(s)) \\[0.75ex]
	\hphantom{d\:\!p(s\:\!;t) = - \, \mathopen{\bigg{\lbrace}} \:\!} + \bar{X}(s)\frac{\partial}{\partial x}\Pi_2(s,\bar{X}(s))\big{\rbrack} \:\!\! \mathclose{\bigg{\rbrace}}\:\!ds \\[0.5ex]
	\hphantom{d\:\!p(s\:\!;t) = \:\!} + q(s\:\!;t)\:\!dW(s), \\[1ex]
p(T\:\!;t) = -\:\!\hat{\hslash}(T\:\!;t)\:\!\hat{\upsilon}'(\bar{X}(T)),
\end{sistema}
\end{equation} 
(where $s\in\mathopen{[}t,T\mathclose{]}$) with, for any $s\in\mathopen{[}t,T\mathclose{]}$ (and $\probP$-\:\!a.s.),
\begin{equation}
\label{eq:applicationuncon1}
p(s\:\!;t) = -\:\!\hslash(s\:\!;t)\:\!\upsilon'(\Pi_2(s,\bar{X}(s))\bar{X}(s)),
\end{equation} 
and
\begin{equation}
\label{eq:applicationuncon2}
\bracks[\big]{\mu - \sigma\:\!\hslash(s\:\!;t)\gamma(s\:\!;t)}p(s\:\!;t) + \sigma\:\!q(s\:\!;t) = 0.
\end{equation} 
If $q(t\:\!;t)$ is a deterministic constant, then $\Py$ is an equilibrium policy having, in particular,
\begin{equation*}
\Pi_2(t,x) = \frac{1}{x} \;\! \Upsilon\big{(}\:\!p(t\:\!;t)\big{)}.
\end{equation*} 
\end{theorem} 

\begin{remark}
\label{rem:remarkapplicationuncon}
Under the conditions~\eqref{eq:applicationuncon1} and~\eqref{eq:applicationuncon2}, the BSDE~\eqref{eq:applicationuncon0} coincides with
\begin{equation}
\label{eq:remarkapplicationuncon}
\begin{sistema}
d\:\!p(s\:\!;t) = -\:\!\braces[\Big]{r\:\!p(s\:\!;t) + \hslash(s\:\!;t)\bracks[\big]{-\:\!\beta(s\:\!;t)p(s\:\!;t) - \gamma(s\:\!;t)q(s\:\!;t)} \:\!\! }\:\!ds + q(s\:\!;t)\:\!dW(s), \\[1ex]
p(T\:\!;t) = -\:\!\hat{\hslash}(T\:\!;t)\:\!\hat{\upsilon}'(\bar{X}(T)),
\end{sistema}
\end{equation} 
(where $s\in\mathopen{[}t,T\mathclose{]}$), from where, thanks again to Remark~\ref{rem:comparison}, we rediscover that
\[
p(\:\!\bm{\cdot}\,;t) < 0
\]
($\probP$-\:\!a.s.), coherently also with~\eqref{eq:applicationuncon1} itself. Furthermore, since $p(t\:\!;t)$ is a deterministic constant (see Remark~\ref{rem:pPdeterministic} or, simply,~\eqref{eq:applicationuncon1}), $q(t\:\!;t)$ is deterministic if and only if $\gamma(t\:\!;t)$ is deterministic (see~\eqref{eq:applicationuncon2}).
\end{remark} 






Now, we introduce the Fenchel\:\!-Legendre transform $\tilde{\upsilon} \colon \mathopen{]}0,\infty\mathclose{[}{_{\:\!y}} \to \mathopen{]}0,\infty\mathclose{[}$ of the (convex) function $- \:\! \upsilon(- \;\! \bm{\cdot} \:\!) \colon \mathopen{]}-\infty,0\mathclose{]} \to \mathopen{]}-\infty,0\mathclose{]}$ given, for any $y \in \mathopen{]}0,\infty\mathclose{[}$, by
\begin{equation}
\label{eq:legendre}
\tilde{\upsilon}(y) \doteq \sup_{z \in \mathopen{[}0,\infty\mathclose{[}} \:\! \braces[\Big]{\! - \:\! y \:\! z + \upsilon(z) \:\!\!} = \upsilon(\:\!\Upsilon(- \:\! y)) - y \:\! \Upsilon(- \:\! y)
\end{equation} 
(see Remark~\ref{rem:Upsilon} and Notation~\ref{not:Upsilon}). Note that $\tilde{\upsilon}$ is (convex and) of class $C^1$\:\!\!, as well as decreasing, with
\[
\tilde{\upsilon}\:\!'(\:\!\bm{\cdot}\:\!) = - \:\! \Upsilon(-  \;\! \bm{\cdot} \:\!).
\]

\begin{assumption}
\label{ass:assumptionapplication}
For any $t \in \mathopen{[}0,T\mathclose{[}$, $\beta(\:\!\bm{\cdot}\,;t)$ and $\gamma(\:\!\bm{\cdot}\,;t)$ are actually deterministic and the functions
\[
\hslash(\:\!\bm{\cdot}\,;t) \gamma(\:\!\bm{\cdot}\,;t), \qquad - \:\! \frac{\hslash'(\:\!\bm{\cdot}\,;t)}{\hslash(\:\!\bm{\cdot}\,;t)} + \hslash(\:\!\bm{\cdot}\,;t) \beta(\:\!\bm{\cdot}\,;t),
\]
defined on $\mathopen{[}t,T\mathclose{]}$, do not depend on $t$ and can be considered as well defined on the whole $\mathopen{[}0,T\mathclose{]}$. Also,
\[
\frac{\hat{\hslash}(T\:\!;t)}{\hslash(T\:\!;t)}
\]
does not depend on $t$.
\end{assumption} 

\begin{remark}
It is not difficult to find non-trivial cases in which Assumption~\ref{ass:assumptionapplication} is satisfied: for instance, by taking, for any $t \in \mathopen{[}0,T\mathclose{[}$, $\beta(\:\!\bm{\cdot}\,;t)$, $\gamma(\:\!\bm{\cdot}\,;t)$ and $\hslash(\:\!\bm{\cdot}\,;t)$ as some exponential functions of the form
\[
e^{\:\! A(\:\!\bm{\cdot}\:\!) + B(t)\:\!\!} ,
\]
where $A \colon \mathopen{[}0,T\mathclose{]}{_s} \to \R$ and $B \colon \mathopen{[}0,T\mathclose{[}{_t} \to \R$ are suitable maps, and then $\hat{\hslash}(\:\!\bm{\cdot}\,;t) = \hslash(\:\!\bm{\cdot}\,;t)$.
\end{remark}

\begin{theorem}
\label{th:hjb}
Suppose Assumption~\ref{ass:assumptionapplication} holds and that there exists a function $V \:\!\! \colon \mathopen{[}0,T\mathclose{]}{_s} \times \mathopen{]}0,\infty\mathclose{[}{_x} \to \R$ of class $C^{1,3}$\:\!\!, with $\frac{\partial \:\! V}{\partial \:\! x} > 0$ and $\frac{\partial^2 \:\! V}{\partial \:\! x^2 \:\!\!} \neq 0$ on $\mathopen{[}0,T\mathclose{]} \times \mathopen{]}0,\infty\mathclose{[}$, which satisfies the PDE
\begin{multline}
\label{eq:hjbeq}
\frac{\partial \:\! V}{\partial \:\! s}(s,x) + r \:\! x \:\! \frac{\partial \:\! V}{\partial \:\! x}(s,x) - \:\! \frac{{\bracks[\big]{\mu - \sigma \:\! \hslash(s\:\!;t) \gamma(s\:\!;t)}}^2 \:\!\!}{2 \:\! \sigma^2\:\!\!} \:\! \frac{\rounds{\frac{\partial \:\! V}{\partial \:\! x}(s,x)}^2 \:\!\!}{\frac{\partial^2 \:\! V}{\partial \:\! x^2 \:\!\!}(s,x)}	\\[0.75ex]
	+ \tilde{\upsilon}{\mathopen{\bigg{(}}}\;\!\!\frac{\partial \:\! V}{\partial \:\! x}(s,x)\!{\mathclose{\bigg{)}}} = \bracks*{- \:\! \frac{\hslash'(s\:\!;t)}{\hslash(s\:\!;t)} + \hslash(s\:\!;t) \beta(s\:\!;t)} V(s,x)
\end{multline} 
(where $s\in\mathopen{[}0,T\mathclose{]}$ and $x\in\mathopen{]}0,\infty\mathclose{[}$) and the boundary condition
\begin{equation}
\label{eq:hjbcond}
V(T,\bm{\cdot}\:\!) = \frac{\hat{\hslash}(T\:\!;t)}{\hslash(T\:\!;t)} \;\! \hat{\upsilon}(\:\!\bm{\cdot}\:\!).
\end{equation} 
Let $\Py = \rounds{\Pi_1,\Pi_2} \colon \mathopen{[}0,T\mathclose{]}{_s} \times \mathopen{]}0,\infty\mathclose{[}{_x} \to \R \times \mathopen{]}0,\infty\mathclose{[}$ be defined, for any $s\in\mathopen{[}0,T\mathclose{]}$ and $x\in\mathopen{]}0,\infty\mathclose{[}$, by
\begin{equation}
\label{eq:hjbPy}
\Pi_1(s,x) \doteq -\:\!\frac{\bracks[\big]{\mu - \sigma\:\!\hslash(s\:\!;t)\gamma(s\:\!;t)} \frac{\partial \:\! V}{\partial \:\! x}(s,x)}{x \;\! \sigma^2 \;\! \frac{\partial^2 \:\! V}{\partial \:\! x^2 \:\!\!}(s,x)}, \qquad \Pi_2(s,x) \doteq \frac{1}{x} \;\! \Upsilon{\mathopen{\bigg{(}}}\!\! - \frac{\partial \:\! V}{\partial \:\! x}(s,x)\!{\mathclose{\bigg{)}}}.
\end{equation} 
If $\Py$ is of class $C^2$\:\!\! w.r.t. $x$ with bounded first and second derivatives even if multiplied by $x$ and if, for any $t\in\mathopen{[}0,T\mathclose{[}$ and $x \in \mathopen{]}0,\infty\mathclose{[}$, $\Py(\:\!\bm{\cdot}\,,\bar{X}(\:\!\bm{\cdot}\:\!)) \in \UUU\bracks{t,T}$ where $\bar{X}(\:\!\bm{\cdot}\:\!)$ is the solution of the FSDE~\eqref{eq:applicationX}, then $\Py$ is an equilibrium policy.
\end{theorem} 
 
\begin{proof}
For computational convenience, we consider $\bm{F} \:\!\! = \rounds{F_1,F_2} \colon \mathopen{[}0,T\mathclose{]}{_s} \times \mathopen{]}0,\infty\mathclose{[}{_x} \to \R \times \mathopen{]}0,\infty\mathclose{[}$ given, for any $s\in\mathopen{[}0,T\mathclose{]}$ and $x\in\mathopen{]}0,\infty\mathclose{[}$, by
\begin{equation}
\label{eq:hjbF}
F_1(s,x) \doteq x \;\! \Pi_1(s,x), \qquad F_2(s,x) \doteq x \;\! \Pi_2(s,x)
\end{equation} 
(see~\eqref{eq:hjbPy}). In particular, for any $s\in\mathopen{[}0,T\mathclose{]}$ and $x\in\mathopen{]}0,\infty\mathclose{[}$,
\[
\frac{\partial \:\! V}{\partial \:\! x}(s,x) = \upsilon'(F_2(s,x)).
\]
For any fixed $t\in\mathopen{[}0,T\mathclose{[}$, we define the process $\rounds{p(\:\!\bm{\cdot}\,;t),q(\:\!\bm{\cdot}\,;t)} \in \LLped{2}{\Omega;\CC\rounds{\mathopen{[}t,T\mathclose{]};\R}}{\F} \times \LLped{2}{t,T;\R}{\F}$ as
\begin{equation}
\label{eq:hjbpq}
p(s\:\!;t) \doteq - \:\! \hslash(s\:\!;t) \frac{\partial \:\! V}{\partial \:\! x}(s,\bar{X}(s)), \qquad q(s\:\!;t) \doteq - \:\! \sigma \:\! \hslash(s\:\!;t) F_1(s,\bar{X}(s)) \frac{\partial^2 \:\! V}{\partial \:\! x^2 \:\!\!}(s,\bar{X}(s))
\end{equation} 
for any $s\in\mathopen{[}t,T\mathclose{]}$. Thereby,
\[
p(T\:\!;t) = -\:\!\hat{\hslash}(T\:\!;t)\:\!\hat{\upsilon}'(\bar{X}(T))
\]
(see~\eqref{eq:hjbcond}) and the conditions~\eqref{eq:applicationuncon1} and~\eqref{eq:applicationuncon2} in Theorem~\ref{th:applicationuncon} are satisfied (see~\eqref{eq:hjbF}). Furthermore, $q(t\:\!;t)$ is a deterministic constant (conditional on the event $\bar{X}(t) = x$). Now the claim is that $\rounds{p(\:\!\bm{\cdot}\,;t),q(\:\!\bm{\cdot}\,;t)}$ is the pair solution of the BSDE~\eqref{eq:remarkapplicationuncon}, which would conclude the proof by virtue of Theorem~\ref{th:applicationuncon} itself (see also Remark~\ref{rem:remarkapplicationuncon}). To this end, we observe that, by~\eqref{eq:applicationX} (and~\eqref{eq:hjbF}),
\[
d\bar{X}(s) = \bracks[\big]{r \:\! \bar{X}(s) + \mu \:\! F_1(s,\bar{X}(s)) - F_2(s,\bar{X}(s))} \:\! ds + \sigma \:\! F_1(s,\bar{X}(s)) \:\! dW(s)
\]
(where $s\in\mathopen{[}0,T\mathclose{]}$) and we apply It\^o's formula to $p(\:\!\bm{\cdot}\,;t)$ (on $\mathopen{[}t,T\mathclose{]}\times\Omega$) obtaining, by~\eqref{eq:hjbpq},
\begin{equation*}
\begin{split}
d\:\!p(s\:\!;t)	&= - \:\! \hslash(s\:\!;t) \bigg{\lbrace} \frac{\hslash'(s\:\!;t)}{\hslash(s\:\!;t)} \frac{\partial \:\! V}{\partial \:\! x}(s,\bar{X}(s)) + \frac{\partial^2 \:\! V}{\partial \:\! s \, \partial \:\! x}(s,\bar{X}(s)) \\[0.5ex]
 			&\hphantom{= - \:\! \hslash(s\:\!;t) \bigg{\lbrace}} + \bracks[\big]{r \:\! \bar{X}(s) + \mu \:\! F_1(s,\bar{X}(s)) - F_2(s,\bar{X}(s))} \frac{\partial^2 \:\! V}{\partial \:\! x^2 \:\!\!}(s,\bar{X}(s)) \\[0.75ex]
 			&\hphantom{= - \:\! \hslash(s\:\!;t) \bigg{\lbrace}} + \frac{\sigma^2 \:\!\!}{2} F_1^{\:\!2}(s,\bar{X}(s)) \frac{\partial^3 \:\! V}{\partial \:\! x^3 \:\!\!}(s,\bar{X}(s)) \:\!\! \bigg{\rbrace} \:\! ds \\[0.5ex]
			&+ q(s\:\!;t) \:\! dW(s)
\end{split}
\end{equation*}
(where $s\in\mathopen{[}t,T\mathclose{]}$). Therefore,~\eqref{eq:remarkapplicationuncon} is fulfilled if and only if, for any $s\in\mathopen{[}t,T\mathclose{]}$ (and $\probP$-\:\!a.s.),
\begin{multline*}
\frac{\partial^2 \:\! V}{\partial \:\! s \, \partial \:\! x}(s,\bar{X}(s)) + r \:\! \bar{X}(s) \frac{\partial^2 \:\! V}{\partial \:\! x^2 \:\!\!}(s,\bar{X}(s)) + r \:\! \frac{\partial \:\! V}{\partial \:\! x}(s,\bar{X}(s)) + \bracks*{\frac{\hslash'(s\:\!;t)}{\hslash(s\:\!;t)} - \hslash(s\:\!;t) \beta(s\:\!;t)} \frac{\partial \:\! V}{\partial \:\! x}(s,\bar{X}(s)) \\[0.5ex]
	+ \braces[\Big]{\bracks[\big]{\mu - \sigma \:\! \hslash(s\:\!;t) \gamma(s\:\!;t)} F_1(s,\bar{X}(s)) - F_2(s,\bar{X}(s))} \frac{\partial^2 \:\! V}{\partial \:\! x^2 \:\!\!}(s,\bar{X}(s)) \\[0.75ex]
		+ \frac{\sigma^2 \:\!\!}{2} F_1^{\:\!2}(s,\bar{X}(s)) \frac{\partial^3 \:\! V}{\partial \:\! x^3 \:\!\!}(s,\bar{X}(s)) = 0,
\end{multline*}
i.e.,
\begin{multline*}
\frac{\partial}{\partial \:\! x} \Bigg{\lbrace} \frac{\partial \:\! V}{\partial \:\! s}(s,x) + r \:\! x \:\! \frac{\partial \:\! V}{\partial \:\! x}(s,x) + \bracks*{\frac{\hslash'(s\:\!;t)}{\hslash(s\:\!;t)} - \hslash(s\:\!;t) \beta(s\:\!;t)} V(s,x)	\\[0.5ex]
	 + \braces[\Big]{\bracks[\big]{\mu - \sigma \:\! \hslash(s\:\!;t) \gamma(s\:\!;t)} F_1(s,x) - F_2(s,x)} \frac{\partial \:\! V}{\partial \:\! x}(s,x) \\[0.75ex]
	 	+ \frac{\sigma^2 \:\!\!}{2} F_1^{\:\!2}(s,x) \frac{\partial^2 \:\! V}{\partial \:\! x^2 \:\!\!}(s,x) + \upsilon(F_2(s,x)) \Bigg{\rbrace}_{ \:\! \bigg{\vert} \:\! x \;\! = \:\! \bar{X}(s) } \! = 0,
\end{multline*}
which holds due to~\eqref{eq:hjbeq} because, also by~\eqref{eq:legendre},
\begin{multline*}
\braces[\Big]{\bracks[\big]{\mu - \sigma \:\! \hslash(s\:\!;t) \gamma(s\:\!;t)} F_1(s,x) - F_2(s,x)} \frac{\partial \:\! V}{\partial \:\! x}(s,x) + \frac{\sigma^2 \:\!\!}{2} F_1^{\:\!2}(s,x) \frac{\partial^2 \:\! V}{\partial \:\! x^2 \:\!\!}(s,x) + \upsilon(F_2(s,x)) \\[1ex]
	= - \:\! \frac{{\bracks[\big]{\mu - \sigma \:\! \hslash(s\:\!;t) \gamma(s\:\!;t)}}^2 \:\!\!}{2 \:\! \sigma^2\:\!\!} \:\! \frac{\rounds{\frac{\partial \:\! V}{\partial \:\! x}(s,x)}^2 \:\!\!}{\frac{\partial^2 \:\! V}{\partial \:\! x^2 \:\!\!}(s,x)} + \tilde{\upsilon}{\mathopen{\bigg{(}}}\;\!\!\frac{\partial \:\! V}{\partial \:\! x}(s,x)\!{\mathclose{\bigg{)}}}
\end{multline*}
for any $s\in\mathopen{[}0,T\mathclose{]}$ and $x\in\mathopen{]}0,\infty\mathclose{[}$.	\qedhere
\end{proof}

The PDE~\eqref{eq:hjbeq} is generally affected by non-uniqueness of the possible solutions: consequently, there could be heterogeneous types of equilibrium policies. Here we take $\lambda \in \mathopen{]}0,1\mathclose{[}$ and, for $x\in\mathopen{[}0,\infty\mathclose{[}$,
\[
\upsilon(x) = \hat{\upsilon}(x) \doteq x^{\lambda}/\lambda
\]
(see  Example~\ref{ex:utilityfunctions}), thus having that, for any $y \in \mathopen{]}0,\infty\mathclose{[}$,
\[
\tilde{\upsilon}(y) = \frac{1 - \lambda}{\lambda} \, y^{- \:\! \frac{\lambda}{1 - \lambda}} \:\!\!
\]
(see~\eqref{eq:legendre}) and therefore, if we look for a solution $V \:\!\! \colon \mathopen{[}0,T\mathclose{]}{_s} \times \mathopen{]}0,\infty\mathclose{[}{_x} \to \R$ of the PDE~\eqref{eq:hjbeq} satisfying the boundary condition~\eqref{eq:hjbcond} which, for any $s\in\mathopen{[}0,T\mathclose{]}$ and $x\in\mathopen{]}0,\infty\mathclose{[}$, can be decomposed as
\[
V(s,x) = \psi(s) \:\! \upsilon(x),
\]
where $\psi \colon \mathopen{[}0,T\mathclose{]}{_s} \to \R$ is a suitable map, then Theorem~\ref{th:hjb} easily reduces to the following result.

\begin{corollary}
\label{cor:hjb}
Suppose Assumption~\ref{ass:assumptionapplication} holds and that there exists a function $\psi \colon \mathopen{[}0,T\mathclose{]}{_s} \to \R$ of class $C^1$\:\!\! which satisfies the ODE
\begin{multline}
\label{eq:corhjb}
\psi'(s) +  \braces*{\lambda \:\! r + \frac{\lambda {\bracks[\big]{\mu - \sigma \:\! \hslash(s\:\!;t) \gamma(s\:\!;t)}}^2 \:\!\!}{2 \:\! (1 - \lambda) \:\! \sigma^2\:\!\!}} \:\! \psi(s)	\\
	+ (1 - \lambda) {\bracks*{\;\! \psi(s) \:\!}}^{\:\! - \:\! \frac{\lambda}{1 - \lambda}} \:\!\! = \bracks*{- \:\! \frac{\hslash'(s\:\!;t)}{\hslash(s\:\!;t)} + \hslash(s\:\!;t) \beta(s\:\!;t)} \psi(s)
\end{multline} 
(where $s\in\mathopen{[}0,T\mathclose{]}$) and the boundary condition
\[
\psi(T) = \frac{\hat{\hslash}(T\:\!;t)}{\hslash(T\:\!;t)}.
\]
Let $\Py = \rounds{\Pi_1,\Pi_2} \colon \mathopen{[}0,T\mathclose{]}{_s} \times \mathopen{]}0,\infty\mathclose{[}{_x} \to \R \times \mathopen{]}0,\infty\mathclose{[}$ be defined, for any $s\in\mathopen{[}0,T\mathclose{]}$ (and $x\in\mathopen{]}0,\infty\mathclose{[}$\:\!), by
\[
\Pi_1(s,x) \doteq \frac{\mu - \sigma\:\!\hslash(s\:\!;t)\gamma(s\:\!;t)}{(1 - \lambda) \:\! \sigma^2 \:\!\!} \:\!, \qquad \Pi_2(s,x) \doteq {\bracks*{\;\! \psi(s) \:\!}}^{\:\! - \:\! \frac{1}{1 - \lambda}} .
\]
Then $\Py$ is an equilibrium policy.
\end{corollary} 

\begin{remark}
Equation~\eqref{eq:hjbeq} in Theorem~\ref{th:hjb} and~\eqref{eq:corhjb} in Corollary~\ref{cor:hjb} are of the same type as those obtained in Theorem~3.2, respectively in Lemma~3.1, in~\cite{ekelandpirvu08}. Hence, existence and uniqueness results of~\eqref{eq:corhjb} can be proved in an analogous, and even easier, way. For more details we refer to~(\cite{ekelandpirvu08}, Proof of Lemma 3.1). Note that the local existence of  $\psi$ close to $T$ is guaranteed; furthermore, we stress that $\psi$ must remain positive.
\end{remark}


\section{Conclusions}
\label{sec:conclusions}


We have formulated the class of time\:\!-inconsistent recursive stochastic optimal control problems, where the notion of optimality is defined by means of subgame\:\!-perfect equilibrium, and we have obtained sufficient and necessary conditions for existence in the form of a (Pontryagin) maximum principle relying on a generalized second-order Hamiltonian function.

Under suitable conditions of analytical-geometric regularity, it is possible to restrict attention to the first-order part of the Hamiltonian alone, as is done with the investment-consumption policies considered in Section~\ref{sec:applicationuncon}.

As regards possible future developments of what has been discussed, the following is highlighted.

Introduce a constraint to the problem and seek at least necessary conditions of existence. Such a constraint could be defined on an expected value---which in turn derives or not from a recursive utility---or it could also be an infinite-dimensional constraint such as
\[
X(T) \in \QQ \subset \Lped{2}{\Omega;\R}{T}
\]
(see, e.g.,~\cite{yong99},~\cite{elkarouipengquenez01} and~\cite{zhuo18}).

We emphasize in this regard that such constrained problems are still quite far from being fully developed: indeed, in the existing literature, stochastic control problems have been studied under the influence of a constraint but, generally, not with recursive utilities. We refer to~\cite{pirvu07} as one of the first notable works in portfolio choice theory with constant-relative-risk aversion (CRRA)-type preferences, for a convex and compact constraint defined through a (pseudo) risk measure such as value at risk ($\VaR$) on a wealth process at a future time instant ``very close'' to the present. Here, the market coefficients are random but independent of the Brownian motion driving the stocks. For a generalization, see~\cite{morenobrombergpirvureveillac13} (CRRA preferences) and also~\cite{huimkellermuller05} and~\cite{cheriditohu11} (martingale methods), among others. 

As far as practical applications are concerned, other portfolio management problems should be explored, with various choices of recursive utility and, possibly, of constraint, moving towards the solution of more general and challenging cases than the one discussed in Section~\ref{sec:applicationuncon}.

A completely new theory could be constructed once the functional $J(\:\!\bm{\cdot}\,;t,x)$ has been modified as
\[
J(\u(\:\!\bm{\cdot}\:\!);t,x) \doteq \E\bigg{\lbrack}\int_t^T\! \ell(s,X(s),\u(s),Y(s\:\!;t),Z(s\:\!;t);t)\,ds + \phi(X(T);t) + \gamma(Y(t\:\!;t);t)\bigg{\rbrack} \:\!\! ,
\]
among others (see, e.g.,~\cite{jizhou06}).

Finally, extensions to the infinite horizon case
\[
T = \infty,
\]
or to random horizons $\tau(\:\!\bm{\cdot}\:\!)$ (stopping times), should be investigated as well.


\appendix


\section{A proof of Lemma~\ref{lemma:liminf}}
\label{subsec:lemmaliminf}


We start from the following notational convention, borrowed from~\cite{hu17}, with regard to which there should be no misunderstandings in this section.

\begin{notation}
\label{not:Hu17}
For $\epsilon\downarrow0$ and $\Lambda^\epsilon(\:\!\bm{\cdot}\,;t) = (\Lambda^\epsilon(s\:\!;t))_{s\in\mathopen{[}t,T\mathclose{]}} \in \LLped{2}{t,T;\R}{\F}$ (possibly a random variable), if
\[
\E\:\!\!{\mathopen{\bigg{(}}\int_t^T\!\abs*{\Lambda^\epsilon(s\:\!;t)}\:\!ds\mathclose{\bigg{)}}}^{\!\!2k}\! = o_{\epsilon\downarrow0}{\big(}\epsilon^{2k}{\big)}
\]
for any $k\in\mathopen{[}1,\infty\mathclose{[}$, then we simply write
\[
\Lambda^\epsilon(\:\!\bm{\cdot}\,;t) = o_{\epsilon\downarrow0}{\big(}\epsilon{\big)}.
\]
\end{notation} 

Fix $\u(\:\!\bm{\cdot}\:\!)\in\UUU\bracks{t,T}$. For $\epsilon\in\mathopen{]}0,T-t\mathclose{[}$ and $s\in\mathopen{[}t,T\mathclose{]}$, let $E_t^{\:\!\epsilon} \equiv E_{t,s}^{\:\!\epsilon}$\:\!\! be as in~\eqref{eq:Etseps} (see Lemma~\ref{lemma:liminf}) and let $\bar{\u}^\epsilon(\:\!\bm{\cdot}\:\!)$ be the spike variation of $\bar{\u}(\:\!\bm{\cdot}\:\!)$ w.r.t. $\u(\:\!\bm{\cdot}\:\!)$ and $E_t^{\:\!\epsilon}$\:\!\!.

Consider the usual approximate variational systems/processes of first and second order of the state process $\bar{X}(\:\!\bm{\cdot}\:\!)$ w.r.t. the control perturbation $\bar{\u}^\epsilon(\:\!\bm{\cdot}\:\!)$: i.e., respectively,
\[
\begin{sistema}
dX_1^\epsilon(s\:\!;t) = b_x(s)X_1^\epsilon(s\:\!;t)\:\!ds + \bracks[\big]{\sigma_x(s)X_1^\epsilon(s\:\!;t) + \delta\sigma(s)\1_{\:\!\!E_t^{\:\!\epsilon}}(s)}dW(s), \\[0.75ex]
X_1^\epsilon(t\:\!;t) = 0,
\end{sistema}
\]
and
\[
\begin{sistema}
dX_2^\epsilon(s\:\!;t) = \bracks*{\:\!b_x(s)X_2^\epsilon(s\:\!;t) + \delta b(s)\1_{\:\!\!E_t^{\:\!\epsilon}}(s) + \frac{1}{2}\:\!b_{xx}(s){X_1^\epsilon(s\:\!;t)}^2\:\!}ds \\[1ex]
	\hphantom{d X_2^\epsilon(s\:\!;t) = } + \big{\lbrack}\:\!\sigma_x(s)X_2^\epsilon(s\:\!;t) + \delta\sigma_x(s)X_1^\epsilon(s\:\!;t)\1_{\:\!\!E_t^{\:\!\epsilon}}(s) \\[1ex]
		\hphantom{d X_2^\epsilon(s\:\!;t) = + \big{\lbrack}\:\!} + \frac{1}{2}\:\!\sigma_{xx}(s){X_1^\epsilon(s\:\!;t)}^2\:\!\big{\rbrack}dW(s), \\[1ex]
X_2^\epsilon(t\:\!;t) = 0,
\end{sistema}
\]
(where $s\in\mathopen{[}t,T\mathclose{]}$). The following result is Theorem~4.4 in~\cite[Chap. 3, Sect. 4]{yongzhou99}, about ``Taylor expansions.'' 

\begin{lemma}
\label{lemma:4.4}
For any $k\in\mathopen{[}1,\infty\mathclose{[}$,
\begin{gather*}
\sup_{s\in\mathopen{[}\:\!t,T\:\!\mathclose{]}}\:\!\! \E\bracks*{{\abs*{X^\epsilon(s) - \bar{X}(s)}}^{2k}} = O{\big(}\epsilon^k{\big)}, \\[0.5ex]
\sup_{s\in\mathopen{[}\:\!t,T\:\!\mathclose{]}}\:\!\! \E\bracks*{{\abs*{X_1^\epsilon(s\:\!;t)}}^{2k}} = O{\big(}\epsilon^k{\big)}, \\[0.5ex]
\sup_{s\in\mathopen{[}\:\!t,T\:\!\mathclose{]}}\:\!\! \E\bracks*{{\abs*{X^\epsilon(s) - \bar{X}(s) - X_1^\epsilon(s\:\!;t)}}^{2k}} = O{\big(}\epsilon^{2k}{\big)}, \\[0.5ex]
\sup_{s\in\mathopen{[}\:\!t,T\:\!\mathclose{]}}\:\!\! \E\bracks*{{\abs*{X_2^\epsilon(s\:\!;t)}}^{2k}} = O{\big(}\epsilon^{2k}{\big)}, \\[0.5ex]
\sup_{s\in\mathopen{[}\:\!t,T\:\!\mathclose{]}}\:\!\! \E\bracks*{{\abs*{X^\epsilon(s) - \bar{X}(s) - X_1^\epsilon(s\:\!;t) - X_2^\epsilon(s\:\!;t)}}^{2k}} = o_{\epsilon\downarrow0}{\big(}\epsilon^{2k}{\big)}.
\end{gather*}
Furthermore,
\begin{multline}
\label{eq:lemma4.4}
h(X^\epsilon(T);t) - h(\bar{X}(T);t) - h_x(\bar{X}(T);t)\bracks[\big]{X_1^\epsilon(T\:\!;t) + X_2^\epsilon(T\:\!;t)} \\[0.75ex]
	- \textstyle{\frac{1}{2}}\:\!h_{xx}(\bar{X}(T);t){X_1^\epsilon(T\:\!;t)}^2\:\!\! = o_{\epsilon\downarrow0}{\big(}\epsilon{\big)}
\end{multline} 
(see Notation~\ref{not:XYZ} and Notation~\ref{not:Hu17}).
\end{lemma} 

\begin{remark}
\label{rem:lemma4.4}
Regarding Lemma~\ref{lemma:4.4}, we point out that~\eqref{eq:lemma4.4} can be rewritten as
\[
Y^\epsilon(T\:\!;t) - \bar{Y}(T\:\!;t) = o_{\epsilon\downarrow0}{\big(}\epsilon{\big)} + p(T\:\!;t)\bracks[\big]{X_1^\epsilon(T\:\!;t) + X_2^\epsilon(T\:\!;t)} + \textstyle{\frac{1}{2}}P(T\:\!;t){X_1^\epsilon(T\:\!;t)}^2\:\!\!
\]
(see Definitions~\ref{def:recutility},~\ref{def:p} and~\ref{def:P}).
\end{remark} 

What we need to properly estimate is, for $\epsilon\downarrow0$,
\[
Y^\epsilon(t\:\!;t) - \bar{Y}(t\:\!;t)
\]
(see~\eqref{eq:liminfY} in Remark~\ref{rem:liminfY}), while, in the sense of Lemma~\ref{lemma:4.4}, we know something useful only about
\[
Y^\epsilon(T\:\!;t) - \bar{Y}(T\:\!;t)
\]
(see Remark~\ref{rem:lemma4.4}). Therefore, the idea is to reconstruct information by going back from this term by means of appropriate BSDEs (using the adjoint processes as in Definitions~\ref{def:p} and~\ref{def:P}).

Grouping w.r.t. $\1_{\:\!\!E_t^{\:\!\epsilon}}(\:\!\bm{\cdot}\:\!)$, $X_1^\epsilon(\:\!\bm{\cdot}\,;t)$, $X_2^\epsilon(\:\!\bm{\cdot}\,;t)$, ${X_1^\epsilon(\:\!\bm{\cdot}\,;t)}^2$\:\!\! and $X_1^\epsilon(\:\!\bm{\cdot}\,;t)\1_{\:\!\!E_t^{\:\!\epsilon}}(\:\!\bm{\cdot}\:\!)$ (on $\mathopen{[}t,T\mathclose{]}\times\Omega$), we have
\begin{equation*}
\begin{split}
 & d\mathopen{\Big{(}}p(\:\!\bm{\cdot}\,;t)\bracks[\big]{X_1^\epsilon(\:\!\bm{\cdot}\,;t) + X_2^\epsilon(\:\!\bm{\cdot}\,;t)} + \textstyle{\frac{1}{2}}P(\:\!\bm{\cdot}\,;t){X_1^\epsilon(\:\!\bm{\cdot}\,;t)}^2\mathclose{\Big{)}}(s)  \\
 &= \bigg{\lbrace} \bracks*{p(s\:\!;t)\delta b(s) + q(s\:\!;t)\delta\sigma(s) + \textstyle{\frac{1}{2}}P(s\:\!;t){\rounds{\delta\sigma(s)}}^2}\:\!\!\1_{\:\!\!E_t^{\:\!\epsilon}}(s) \\
 &\hphantom{= \bigg{\lbrace}} + \bracks[\big]{p(s\:\!;t)b_x(s) + q(s\:\!;t)\sigma_x(s) - g(s,p(s\:\!;t),q(s\:\!;t);t)}\:\!\!\bracks[\big]{X_1^\epsilon(s\:\!;t) + X_2^\epsilon(s\:\!;t)} \\[1.125ex]
 &\hphantom{= \bigg{\lbrace}} + \textstyle{\frac{1}{2}}\Big{\lbrack}p(s\:\!;t)b_{xx}(s) + q(s\:\!;t)\sigma_{xx}(s) + 2P(s\:\!;t)b_x(s) + 2Q(s\:\!;t)\sigma_x(s) \\
 &\hphantom{= \bigg{\lbrace} + \textstyle{\frac{1}{2}} \Big{\lbrack} \:\!\!\!} + P(s\:\!;t){\sigma_x(s)}^2\:\! - G(s,P(s\:\!;t),Q(s\:\!;t);t)\Big{\rbrack}{X_1^\epsilon(s\:\!;t)}^2 \\
 &\hphantom{= \bigg{\lbrace}} + \bracks[\big]{q(s\:\!;t)\delta\sigma_x(s) + P(s\:\!;t)\sigma_x(s)\delta\sigma(s) + Q(s\:\!;t)\delta\sigma(s)}X_1^\epsilon(s\:\!;t)\1_{\:\!\!E_t^{\:\!\epsilon}}(s) \bigg{\rbrace}\:\!ds \\
 &+ \bigg{\lbrace} p(s\:\!;t)\delta\sigma(s)\1_{\:\!\!E_t^{\:\!\epsilon}}(s) \\
 &\hphantom{+ \bigg{\lbrace}} + \bracks[\big]{p(s\:\!;t)\sigma_x(s) + q(s\:\!;t)}\:\!\!\bracks[\big]{X_1^\epsilon(s\:\!;t) + X_2^\epsilon(s\:\!;t)} \\[1.125ex]
 &\hphantom{+ \bigg{\lbrace}} + \textstyle{\frac{1}{2}}\bracks[\big]{p(s\:\!;t)\sigma_{xx}(s) + 2P(s\:\!;t)\sigma_x(s) + Q(s\:\!;t)}{X_1^\epsilon(s\:\!;t)}^2 \\
 &\hphantom{+ \bigg{\lbrace}} + \bracks[\big]{p(s\:\!;t)\delta\sigma_x(s) + P(s\:\!;t)\delta\sigma(s)}X_1^\epsilon(s\:\!;t)\1_{\:\!\!E_t^{\:\!\epsilon}}(s) \bigg{\rbrace}\:\!dW(s)
\end{split}
\end{equation*}
(see also Notation~\ref{not:phi}). In particular, the maps $g$ and $G$ do not appear in the $W$\:\!\!-\:\!term.	

For the sake of brevity, we denote, for any $s\in\mathopen{[}t,T\mathclose{]}$,
\begin{gather*}
A_1(s\:\!;t) \eqnot p(s\:\!;t)\delta b(s) + q(s\:\!;t)\delta\sigma(s) + \textstyle{\frac{1}{2}}P(s\:\!;t){\rounds{\delta\sigma(s)}}^2\:\!\!, \\[0.25ex]
A_2(s\:\!;t) \eqnot p(s\:\!;t)b_x(s) + q(s\:\!;t)\sigma_x(s) - g(s,p(s\:\!;t),q(s\:\!;t);t), \\[0.25ex]
A_3(s\:\!;t) \eqnot p(s\:\!;t)b_{xx}(s) + q(s\:\!;t)\sigma_{xx}(s) + 2P(s\:\!;t)b_x(s) + 2Q(s\:\!;t)\sigma_x(s) \\
	+ P(s\:\!;t){\sigma_x(s)}^2\:\!\! - G(s,P(s\:\!;t),Q(s\:\!;t);t), \\[0.25ex]
A_4(s\:\!;t) \eqnot q(s\:\!;t)\delta\sigma_x(s) + P(s\:\!;t)\sigma_x(s)\delta\sigma(s) + Q(s\:\!;t)\delta\sigma(s), \\[0.25ex]
\Theta_1(s\:\!;t) \eqnot p(s\:\!;t)\sigma_x(s) + q(s\:\!;t) ,\\[0.25ex]
\Theta_2(s\:\!;t) \eqnot p(s\:\!;t)\sigma_{xx}(s) + 2P(s\:\!;t)\sigma_x(s) + Q(s\:\!;t), \\[0.25ex]
\Theta_3(s\:\!;t) \eqnot p(s\:\!;t)\delta\sigma_x(s) + P(s\:\!;t)\delta\sigma(s)
\end{gather*}
($p(s\:\!;t)\delta\sigma(s)$ is already short enough). We  focus here on the coefficients of $\1_{\:\!\!E_t^{\:\!\epsilon}}(\:\!\bm{\cdot}\:\!)$.

\begin{remark}
$A_4(\:\!\bm{\cdot}\,;t)X_1^\epsilon(\:\!\bm{\cdot}\,;t)\1_{\:\!\!E_t^{\:\!\epsilon}}(\:\!\bm{\cdot}\:\!) = o_{\epsilon\downarrow0}(\epsilon)$, which does not hold true for $\Theta_3(\:\!\bm{\cdot}\,;t)$ (see~\cite{hu17}).
\end{remark}

Putting it all together, for any $s\in\mathopen{[}t,T\mathclose{]}$,
\begin{equation*}
\begin{split}
Y^\epsilon(s\:\!;t) 
	&= h(\bar{X}(T);t) + o_{\epsilon\downarrow0}(\epsilon) + p(s\:\!;t)\bracks[\big]{X_1^\epsilon(s\:\!;t) + X_2^\epsilon(s\:\!;t)} + \textstyle{\frac{1}{2}}P(s\:\!;t){X_1^\epsilon(s\:\!;t)}^2 \\[1.125ex]
	&+ \int_{\:\!\!s}^T\! \Big{\lbrack}f(r,X^\epsilon(r),\u^\epsilon(r),Y^\epsilon(r\:\!;t),Z^\epsilon(r\:\!;t);t) + A_1(r\:\!;t)\1_{\:\!\!E_t^{\:\!\epsilon}}(r) \\
	&\hphantom{+ \int_{\:\!\!s}^T\! \Big{\lbrack}} + A_2(r\:\!;t)\bracks[\big]{X_1^\epsilon(r\:\!;t) + X_2^\epsilon(r\:\!;t)} + \textstyle{\frac{1}{2}}A_3(r\:\!;t){X_1^\epsilon(r\:\!;t)}^2\Big{\rbrack}\:\!dr \\
	&- \int_{\:\!\!s}^T\! \Big{\lbrack}Z^\epsilon(r\:\!;t) - \Big{(}p(r\:\!;t)\delta\sigma(r)\1_{\:\!\!E_t^{\:\!\epsilon}}(r) + \Theta_1(r\:\!;t)\bracks[\big]{X_1^\epsilon(r\:\!;t) + X_2^\epsilon(r\:\!;t)} \\
	&\hphantom{+ \int_{\:\!\!s}^T\! \Big{\lbrack}} + \textstyle{\frac{1}{2}}\:\!\Theta_2(r\:\!;t){X_1^\epsilon(r\:\!;t)}^2\:\!\! + \Theta_3(r\:\!;t)X_1^\epsilon(r\:\!;t)\1_{\:\!\!E_t^{\:\!\epsilon}}(r)\Big{)}\Big{\rbrack}\:\!dW(r).
\end{split}
\end{equation*}
The key idea is to be able to narrow our analysis to the integral in $dr$. So, if we denote, for $\epsilon\downarrow0$,
\[
\tilde{Y}^\epsilon(\:\!\bm{\cdot}\,;t) \eqnot Y^\epsilon(\:\!\bm{\cdot}\,;t) - \Big{(}p(\:\!\bm{\cdot}\,;t)\bracks[\big]{X_1^\epsilon(\:\!\bm{\cdot}\,;t) + X_2^\epsilon(\:\!\bm{\cdot}\,;t)} + \textstyle{\frac{1}{2}}P(\:\!\bm{\cdot}\,;t){X_1^\epsilon(\:\!\bm{\cdot}\,;t)}^2\Big{)}
\]
and
\begin{equation*}
\begin{split}
\tilde{Z}^\epsilon(\:\!\bm{\cdot}\,;t) &\eqnot Z^\epsilon(\:\!\bm{\cdot}\,;t) - \Big{(}p(\:\!\bm{\cdot}\,;t)\delta\sigma(\:\!\bm{\cdot}\:\!)\1_{\:\!\!E_t^{\:\!\epsilon}}(\:\!\bm{\cdot}\:\!) + \Theta_1(\:\!\bm{\cdot}\,;t)\bracks[\big]{X_1^\epsilon(\:\!\bm{\cdot}\,;t) + X_2^\epsilon(\:\!\bm{\cdot}\,;t)} \\
	&\hphantom{\eqnot \;} + \textstyle{\frac{1}{2}}\:\!\Theta_2(\:\!\bm{\cdot}\,;t){X_1^\epsilon(\:\!\bm{\cdot}\,;t)}^2 + \Theta_3(\:\!\bm{\cdot}\,;t)X_1^\epsilon(\:\!\bm{\cdot}\,;t)\1_{\:\!\!E_t^{\:\!\epsilon}}(\:\!\bm{\cdot}\:\!)\Big{)}
\end{split}
\end{equation*}
(on $\mathopen{[}t,T\mathclose{]}\times\Omega$), then we can write, for any $s\in\mathopen{[}t,T\mathclose{]}$,
\begin{equation*}
\begin{split}
\tilde{Y}^\epsilon(s\:\!;t) &= h(\bar{X}(T);t) + o_{\epsilon\downarrow0}(\epsilon) \\[1.125ex]
 &+ \int_{\:\!\!s}^T\! \Big{\lbrack}f(r,X^\epsilon(r),\u^\epsilon(r),Y^\epsilon(r\:\!;t),Z^\epsilon(r\:\!;t);t) + A_1(r\:\!;t)\1_{\:\!\!E_t^{\:\!\epsilon}}(r) \\
 &\hphantom{+ \int_{\:\!\!s}^T\! \Big{\lbrack}} + A_2(r\:\!;t)\bracks[\big]{X_1^\epsilon(r\:\!;t) + X_2^\epsilon(r\:\!;t)} + \textstyle{\frac{1}{2}}A_3(r\:\!;t){X_1^\epsilon(r\:\!;t)}^2\Big{\rbrack}\:\!dr \\
 &- \int_{\:\!\!s}^T\! \tilde{Z}^\epsilon(r\:\!;t)\,dW(r).
\end{split}
\end{equation*}
For $\epsilon\downarrow0$, we define an approximate variational pair process of first order of $\rounds{\bar{Y}(\:\!\bm{\cdot}\,;t),\bar{Z}(\:\!\bm{\cdot}\,;t)}$ as
\[
Y_1^\epsilon(\:\!\bm{\cdot}\,;t) \eqnot \tilde{Y}^\epsilon(\:\!\bm{\cdot}\,;t) - \bar{Y}(\:\!\bm{\cdot}\,;t)
\]
and
\[
Z_1^\epsilon(\:\!\bm{\cdot}\,;t) \eqnot \tilde{Z}^\epsilon(\:\!\bm{\cdot}\,;t) - \bar{Z}(\:\!\bm{\cdot}\,;t)
\]
(on $\mathopen{[}t,T\mathclose{]}\times\Omega$), and we then get  that, for any $s\in\mathopen{[}t,T\mathclose{]}$,
\begin{equation*}
\begin{split}
Y_1^\epsilon(s\:\!;t) &= o_{\epsilon\downarrow0}(\epsilon) \\[1.125ex]
	&+ \int_{\:\!\!s}^T\! \Big{\lbrack} f(r,X^\epsilon(r),\u^\epsilon(r),Y^\epsilon(r\:\!;t),Z^\epsilon(r\:\!;t);t) \\[1.125ex]
	&\hphantom{+ \int_{\:\!\!s}^T\! \Big{\lbrack}} - f(r,\bar{X}(r),\bar{\u}(r),\bar{Y}(r\:\!;t),\bar{Z}(r\:\!;t);t) + A_1(r\:\!;t)\1_{\:\!\!E_t^{\:\!\epsilon}}(r) \\[1.125ex]
	&\hphantom{+ \int_{\:\!\!s}^T\! \Big{\lbrack}} + A_2(r\:\!;t)\bracks[\big]{X_1^\epsilon(r\:\!;t) + X_2^\epsilon(r\:\!;t)} + \textstyle{\frac{1}{2}}A_3(r\:\!;t){X_1^\epsilon(r\:\!;t)}^2\Big{\rbrack}\:\!dr \\[0.25ex]
	&- \int_{\:\!\!s}^T\! Z_1^\epsilon(r\:\!;t)\,dW(r).
\end{split}
\end{equation*}


Finally, for $\epsilon\downarrow0$, let $\rounds{\Delta Y^\epsilon(\:\!\bm{\cdot}\,;t),\Delta\:\!Z^\epsilon(\:\!\bm{\cdot}\,;t)}$ be such that $\Delta Y^\epsilon(T\:\!;t) = 0$, and, for any $s\in\mathopen{[}t,T\mathclose{]}$,
\begin{equation*}
\begin{split}
\Delta Y^\epsilon(s\:\!;t) &= \int_{\:\!\!s}^T\! \bigg{\lbrace}f_y(r\:\!;t)\Delta Y^\epsilon(r\:\!;t) + f_z(r\:\!;t)\Delta\:\!Z^\epsilon(r\:\!;t) \\[1.125ex]
 &\hphantom{= \int_{\:\!\!s}^T\! \bigg{\lbrace}} + \Big{\lbrack}p(r\:\!;t)\delta b(r) + q(r\:\!;t)\delta\sigma(r) + \textstyle{\frac{1}{2}}P(r\:\!;t){\rounds{\delta\sigma(r)}}^2\:\!\! \\[1.5ex]
 &\hphantom{= \int_{\:\!\!s}^T\! \bigg{\lbrace} + \Big{\lbrack}} + f(r,\bar{X}(r),\u(r),\bar{Y}(r\:\!;t),\bar{Z}(r\:\!;t) + p(r\:\!;t)\delta\sigma(r);t) \\[1.125ex]
 &\hphantom{= \int_{\:\!\!s}^T\! \bigg{\lbrace} + \Big{\lbrack}} - f(r,\bar{X}(r),\bar{\u}(r),\bar{Y}(r\:\!;t),\bar{Z}(r\:\!;t);t)\Big{\rbrack}\:\!\!\1_{\:\!\!E_t^{\:\!\epsilon}}(r)\bigg{\rbrace}\:\!dr \\[0.25ex]
 &- \int_{\:\!\!s}^T\! \Delta\:\!Z^\epsilon(r\:\!;t)\,dW(r)
\end{split}
\end{equation*}
(again, $\Delta Y^\epsilon(t\:\!;t)$ is deterministic). The following result is essentially Theorem 1 in~\cite{hu17}.
\begin{lemma}[Hu (2017)]
\label{lemma:hu17}
W.r.t. Notation~\ref{not:Hu17},
\begin{gather*}
Y_1^\epsilon(\:\!\bm{\cdot}\,;t) - \Delta Y^\epsilon(\:\!\bm{\cdot}\,;t) = o_{\epsilon\downarrow0}(\epsilon), \\[0.75ex]
Z_1^\epsilon(\:\!\bm{\cdot}\,;t) - \Delta\:\!Z^\epsilon(\:\!\bm{\cdot}\,;t) = o_{\epsilon\downarrow0}(\epsilon).
\end{gather*}
\end{lemma} 

As a remarkable corollary of Lemma~\ref{lemma:hu17}, we get
\begin{equation}
\label{eq:Yeps}
Y^\epsilon(t\:\!;t) - \bar{Y}(t\:\!;t) = \Delta Y^\epsilon(t\:\!;t) + o_{\epsilon\downarrow0}(\epsilon)
\end{equation} 
(as well as $Z^\epsilon(t\:\!;t) = \bar{Z}(t\:\!;t) + \Delta\:\!Z^\epsilon(t\:\!;t) + p(t\:\!;t)\delta\sigma(t) + o_{\epsilon\downarrow0}(\epsilon)$).

Furthermore, if $\kappa(\:\!\bm{\cdot}\,;t)$ is the change of numéraire defined through~\eqref{eq:kappa} of Definition~\ref{def:kappa}, then	
\begin{multline}
\label{eq:dYeps}
\Delta Y^\epsilon(t\:\!;t) = \E\!\bigg{\lbrack}\!\int_t^T\! \kappa(r\:\!;t)\Big{\lbrack}p(r\:\!;t)\delta b(r) + q(r\:\!;t)\delta\sigma(r) + \textstyle{\frac{1}{2}}P(r\:\!;t){\rounds{\delta\sigma(r)}}^2\:\!\! \\[1ex]
 	\qquad\qquad + f(r,\bar{X}(r),\u(r),\bar{Y}(r\:\!;t),\bar{Z}(r\:\!;t) + p(r\:\!;t)\delta\sigma(r);t) \\[0.75ex]
		- f(r,\bar{X}(r),\bar{\u}(r),\bar{Y}(r\:\!;t),\bar{Z}(r\:\!;t);t)\Big{\rbrack}\:\!\!\1_{\:\!\!E_t^{\:\!\epsilon}}(r)\,dr\bigg{\rbrack}.
\end{multline} 
In conclusion, for any $s\in\mathopen{[}t,T\mathclose{]}$, by~\eqref{eq:Yeps},~\eqref{eq:dYeps} and the classic Lebesgue differentiation theorem,
\begin{equation*}
\begin{split}
\liminf_{\epsilon\downarrow0}\:\!\frac{J(\bar{\u}^\epsilon(\:\!\bm{\cdot}\:\!);t,x) - J(\bar{\u}(\:\!\bm{\cdot}\:\!);t,x)}{\epsilon} &= \lim_{\epsilon\downarrow0}\:\!\frac{\Delta Y^\epsilon(t\:\!;t)}{\epsilon} \\[1.25ex]
 &= \E\bracks[\big]{\kappa(s\:\!;t)\:\!\delta\HH(s\:\!;t,\u(s))}
\end{split}
\end{equation*}
(see Definitions~\ref{def:H} and~\ref{def:HH} and Notation~\ref{not:HH}), as it was our aim to prove.

\begin{remark}
$\kappa(t\:\!;t)$ could have been defined, in~\eqref{eq:kappa}, as a constant $> 0$ not necessarily equal to $1$.
\end{remark}


\begin{description}[leftmargin=*]
	
	\item[Acknowledgements.] The authors are immensely grateful to Daniele Cassani for his flawless support and to Ivar Ekeland and Santiago Moreno\:\!-Bromberg for several inspiring discussions. They also wish to thank the anonymous Referees for their comments which contributed to improve the paper.

\end{description}



\bibliographystyle{elsarticle-harv} 
\bibliography{elsarticle-template-harv}

@article {bjorkmurgoci14,
    AUTHOR = {Bj\"{o}rk, Tomas and Murgoci, Agatha},
     TITLE = {A theory of {M}arkovian time\:\!-inconsistent stochastic control
              in discrete time},
   JOURNAL = {Finance Stoch.},
  FJOURNAL = {Finance and Stochastics},
    VOLUME = {18},
      YEAR = {2014},
    NUMBER = {3},
     PAGES = {545--592},
      ISSN = {0949-2984},
   MRCLASS = {49L20 (60J05 60J20 91G10 91G80)},
  MRNUMBER = {3232016},
       DOI = {10.1007/s00780-014-0234-y},
       URL = {https://doi.org/10.1007/s00780-014-0234-y},
}

@article {bjorkmurgocizhou14,
    AUTHOR = {Bj\"{o}rk, Tomas and Murgoci, Agatha and Zhou, Xun Yu},
     TITLE = {Mean-variance portfolio optimization with state-dependent risk
              aversion},
   JOURNAL = {Math. Finance},
  FJOURNAL = {Mathematical Finance. An International Journal of Mathematics,
              Statistics and Financial Economics},
    VOLUME = {24},
      YEAR = {2014},
    NUMBER = {1},
     PAGES = {1--24},
      ISSN = {0960-1627},
   MRCLASS = {91G10},
  MRNUMBER = {3157686},
       DOI = {10.1111/j.1467-9965.2011.00515.x},
       URL = {https://doi.org/10.1111/j.1467-9965.2011.00515.x},
}

@article {HaSIAM21,
    AUTHOR = {Hamaguchi, Yushi},
     TITLE = {Time-inconsistent consumption-investment problems in
              incomplete markets under general discount functions},
   JOURNAL = {SIAM J. Control Optim.},
  FJOURNAL = {SIAM Journal on Control and Optimization},
    VOLUME = {59},
      YEAR = {2021},
    NUMBER = {3},
     PAGES = {2121--2146},
      ISSN = {0363-0129},
   MRCLASS = {91G10 (49J55 60H10 93E20)},
  MRNUMBER = {4272906},
MRREVIEWER = {Xiaomin Shi},
       DOI = {10.1137/19M1303782},
       URL = {https://doi.org/10.1137/19M1303782},
}

@article {bjorkkhapkomurgoci17,
    AUTHOR = {Bj\"{o}rk, Tomas and Khapko, Mariana and Murgoci, Agatha},
     TITLE = {On time\:\!-inconsistent stochastic control in continuous time},
   JOURNAL = {Finance Stoch.},
  FJOURNAL = {Finance and Stochastics},
    VOLUME = {21},
      YEAR = {2017},
    NUMBER = {2},
     PAGES = {331--360},
      ISSN = {0949-2984},
   MRCLASS = {49L20 (49N90 60J70 91A10 91A80 91B02 91B25)},
  MRNUMBER = {3626618},
       DOI = {10.1007/s00780-017-0327-5},
       URL = {https://doi.org/10.1007/s00780-017-0327-5},
}

@article {brianddelyonhupardouxstoica03,
    AUTHOR = {Briand, Ph. and Delyon, B. and Hu, Y. and Pardoux, E. and
              Stoica, L.},
     TITLE = {{$L^p$}\:\!\! solutions of backward stochastic differential
              equations},
   JOURNAL = {Stochastic Process. Appl.},
  FJOURNAL = {Stochastic Processes and their Applications},
    VOLUME = {108},
      YEAR = {2003},
    NUMBER = {1},
     PAGES = {109--129},
      ISSN = {0304-4149},
   MRCLASS = {60H10},
  MRNUMBER = {2008603},
MRREVIEWER = {Olivier Raimond},
       DOI = {10.1016/S0304-4149(03)00089-9},
       URL = {https://doi.org/10.1016/S0304-4149(03)00089-9},
}

@article {cheriditohu11,
    AUTHOR = {Cheridito, Patrick and Hu, Ying},
     TITLE = {Optimal consumption and investment in incomplete markets with
              general constraints},
   JOURNAL = {Stoch. Dyn.},
  FJOURNAL = {Stochastics and Dynamics},
    VOLUME = {11},
      YEAR = {2011},
    NUMBER = {2-3},
     PAGES = {283--299},
      ISSN = {0219-4937},
   MRCLASS = {91G10 (60H10 93E20)},
  MRNUMBER = {2836526},
MRREVIEWER = {Tomoyuki Ichiba},
       DOI = {10.1142/S0219493711003280},
       URL = {https://doi.org/10.1142/S0219493711003280},
}

@article {debusschefuhrmantessitore07,
    AUTHOR = {Debussche, Arnaud and Fuhrman, Marco and Tessitore, Gianmario},
     TITLE = {Optimal control of a stochastic heat equation with
              boundary-noise and boundary-control},
   JOURNAL = {ESAIM Control Optim. Calc. Var.},
  FJOURNAL = {ESAIM. Control, Optimisation and Calculus of Variations},
    VOLUME = {13},
      YEAR = {2007},
    NUMBER = {1},
     PAGES = {178--205},
      ISSN = {1292-8119},
   MRCLASS = {49L20 (35K20 35R60 60H15 60H30 93E20)},
  MRNUMBER = {2282108},
MRREVIEWER = {Nicolae H. Pavel},
       DOI = {10.1051/cocv:2007001},
       URL = {https://doi.org/10.1051/cocv:2007001},
}

@article {duffieepstein92,
    AUTHOR = {Duffie, Darrell and Epstein, Larry G.},
     TITLE = {Stochastic differential utility},
      NOTE = {With an appendix by the authors and C. Skiadas},
   JOURNAL = {Econometrica},
  FJOURNAL = {Econometrica. Journal of the Econometric Society},
    VOLUME = {60},
      YEAR = {1992},
    NUMBER = {2},
     PAGES = {353--394},
      ISSN = {0012-9682},
   MRCLASS = {90A43 (49L20 90A10)},
  MRNUMBER = {1162620},
       DOI = {10.2307/2951600},
       URL = {https://doi.org/10.2307/2951600},
}

@Misc{ekelandlazrak06,
  title =	 {Being serious about non\:\!-commitment: subgame perfect equilibrium in continuous time},
  author =	 {Ekeland, Ivar and Lazrak, Ali},
  url =	 {http://arxiv.org/abs/math/0604264},
  year =	 {2006}
}

@article {ekelandmbodjipirvu12,
    AUTHOR = {Ekeland, Ivar and Mbodji, Oumar and Pirvu, Traian A.},
     TITLE = {Time\:\!-consistent portfolio management},
   JOURNAL = {SIAM J. Financial Math.},
  FJOURNAL = {SIAM Journal on Financial Mathematics},
    VOLUME = {3},
      YEAR = {2012},
    NUMBER = {1},
     PAGES = {1--32},
   MRCLASS = {91G10 (60H20 60H30 91B30)},
  MRNUMBER = {2968026},
MRREVIEWER = {George Stoica},
       DOI = {10.1137/100810034},
       URL = {https://doi.org/10.1137/100810034},
}

@article {ekelandpirvu08,
    AUTHOR = {Ekeland, Ivar and Pirvu, Traian A.},
     TITLE = {Investment and consumption without commitment},
   JOURNAL = {Math. Financ. Econ.},
  FJOURNAL = {Mathematics and Financial Economics},
    VOLUME = {2},
      YEAR = {2008},
    NUMBER = {1},
     PAGES = {57--86},
      ISSN = {1862-9679},
   MRCLASS = {91G10 (60G44 60H30)},
  MRNUMBER = {2461340},
MRREVIEWER = {Tommi P. Sottinen},
       DOI = {10.1007/s11579-008-0014-6},
       URL = {https://doi.org/10.1007/s11579-008-0014-6},
}

@article {elkarouipengquenez97,
    AUTHOR = {El Karoui, N. and Peng, S. and Quenez, M. C.},
     TITLE = {Backward stochastic differential equations in finance},
   JOURNAL = {Math. Finance},
  FJOURNAL = {Mathematical Finance. An International Journal of Mathematics,
              Statistics and Financial Economics},
    VOLUME = {7},
      YEAR = {1997},
    NUMBER = {1},
     PAGES = {1--71},
      ISSN = {0960-1627},
   MRCLASS = {90A09 (60H10)},
  MRNUMBER = {1434407},
MRREVIEWER = {R\"{u}diger Kiesel},
       DOI = {10.1111/1467-9965.00022},
       URL = {https://doi.org/10.1111/1467-9965.00022},
}

@article {elkarouipengquenez01,
    AUTHOR = {El Karoui, N. and Peng, S. and Quenez, M. C.},
     TITLE = {A dynamic maximum principle for the optimization of recursive
              utilities under constraints},
   JOURNAL = {Ann. Appl. Probab.},
  FJOURNAL = {The Annals of Applied Probability},
    VOLUME = {11},
      YEAR = {2001},
    NUMBER = {3},
     PAGES = {664--693},
      ISSN = {1050-5164},
   MRCLASS = {49K45 (60J60 91B16 93E20)},
  MRNUMBER = {1865020},
MRREVIEWER = {Lorenzo Peccati},
       DOI = {10.1214/aoap/1015345345},
       URL = {https://doi.org/10.1214/aoap/1015345345},
}

@article {gundelweber08,
    AUTHOR = {Gundel, Anne and Weber, Stefan},
     TITLE = {Utility maximization under a shortfall risk constraint},
   JOURNAL = {J. Math. Econom.},
  FJOURNAL = {Journal of Mathematical Economics},
    VOLUME = {44},
      YEAR = {2008},
    NUMBER = {11},
     PAGES = {1126--1151},
      ISSN = {0304-4068},
   MRCLASS = {91G10 (60G44)},
  MRNUMBER = {2456472},
       DOI = {10.1016/j.jmateco.2008.01.002},
       URL = {https://doi.org/10.1016/j.jmateco.2008.01.002},
}

@article {hu17,
    AUTHOR = {Hu, Mingshang},
     TITLE = {Stochastic global maximum principle for optimization with
              recursive utilities},
   JOURNAL = {Probab. Uncertain. Quant. Risk},
  FJOURNAL = {Probability, Uncertainty and Quantitative Risk},
    VOLUME = {2},
      YEAR = {2017},
     PAGES = {Paper No. 1, 20},
      ISSN = {2095-9672},
   MRCLASS = {49K45 (60H10 93E20)},
  MRNUMBER = {3625730},
       DOI = {10.1186/s41546-017-0014-7},
       URL = {https://doi.org/10.1186/s41546-017-0014-7},
}

@article {huimkellermuller05,
    AUTHOR = {Hu, Ying and Imkeller, Peter and M\"{u}ller, Matthias},
     TITLE = {Utility maximization in incomplete markets},
   JOURNAL = {Ann. Appl. Probab.},
  FJOURNAL = {The Annals of Applied Probability},
    VOLUME = {15},
      YEAR = {2005},
    NUMBER = {3},
     PAGES = {1691--1712},
      ISSN = {1050-5164},
   MRCLASS = {91B28 (60H10 60H30 91B16 91B70 93E20)},
  MRNUMBER = {2152241},
MRREVIEWER = {Vivek S. Borkar},
       DOI = {10.1214/105051605000000188},
       URL = {https://doi.org/10.1214/105051605000000188},
}

@article {hujinzhou12,
    AUTHOR = {Hu, Ying and Jin, Hanqing and Zhou, Xun Yu},
     TITLE = {Time\:\!-inconsistent stochastic linear-quadratic control},
   JOURNAL = {SIAM J. Control Optim.},
  FJOURNAL = {SIAM Journal on Control and Optimization},
    VOLUME = {50},
      YEAR = {2012},
    NUMBER = {3},
     PAGES = {1548--1572},
      ISSN = {0363-0129},
   MRCLASS = {49N10 (60H10 91B02 91G10 93E20)},
  MRNUMBER = {2968066},
MRREVIEWER = {Juan Li},
       DOI = {10.1137/110853960},
       URL = {https://doi.org/10.1137/110853960},
}

@article {hujinzhou17,
    AUTHOR = {Hu, Ying and Jin, Hanqing and Zhou, Xun Yu},
     TITLE = {Time\:\!-inconsistent stochastic linear-quadratic control:
              characterization and uniqueness of equilibrium},
   JOURNAL = {SIAM J. Control Optim.},
  FJOURNAL = {SIAM Journal on Control and Optimization},
    VOLUME = {55},
      YEAR = {2017},
    NUMBER = {2},
     PAGES = {1261--1279},
      ISSN = {0363-0129},
   MRCLASS = {91B51 (60H10 91G10 93E20)},
  MRNUMBER = {3639569},
       DOI = {10.1137/15M1019040},
       URL = {https://doi.org/10.1137/15M1019040},
}

@article {imkellerdosreis10,
    AUTHOR = {Imkeller, Peter and Dos Reis, Gon\c{c}alo},
     TITLE = {Path regularity and explicit convergence rate for {BSDE} with
              truncated quadratic growth},
   JOURNAL = {Stochastic Process. Appl.},
  FJOURNAL = {Stochastic Processes and their Applications},
    VOLUME = {120},
      YEAR = {2010},
    NUMBER = {3},
     PAGES = {348--379},
      ISSN = {0304-4149},
   MRCLASS = {60H07 (60G17 60H10 60H35 65C30)},
  MRNUMBER = {2584898},
MRREVIEWER = {Ali S\"{u}leyman \"{U}st\"{u}nel},
       DOI = {10.1016/j.spa.2009.11.004},
       URL = {https://doi.org/10.1016/j.spa.2009.11.004},
}

@article {jizhou06,
    AUTHOR = {Ji, Shaolin and Zhou, Xun Yu},
     TITLE = {A maximum principle for stochastic optimal control with
              terminal state constraints, and its applications},
   JOURNAL = {Commun. Inf. Syst.},
  FJOURNAL = {Communications in Information and Systems},
    VOLUME = {6},
      YEAR = {2006},
    NUMBER = {4},
     PAGES = {321--337},
      ISSN = {1526-7555},
   MRCLASS = {93E20 (49K45)},
  MRNUMBER = {2346931},
MRREVIEWER = {Brahim Mezerdi},
       URL = {http://projecteuclid.org/euclid.cis/1183729000},
}

@book {mamorelyong99,
    AUTHOR = {Ma, Jin and Yong, Jiongmin},
     TITLE = {Forward-backward stochastic differential equations and their
              applications},
    SERIES = {Lecture Notes in Mathematics},
    VOLUME = {1702},
 PUBLISHER = {Springer\:\!-Verlag, Berlin},
      YEAR = {1999},
     PAGES = {xiv+270},
      ISBN = {3-540-65960-9},
   MRCLASS = {60H10},
  MRNUMBER = {1704232},
MRREVIEWER = {Fabio Antonelli},
}

@article {morenobrombergpirvureveillac13,
    AUTHOR = {Moreno\:\!-Bromberg, Santiago and Pirvu, Traian A. and R\'{e}veillac,
              Anthony},
     TITLE = {C{RRA} utility maximization under dynamic risk constraints},
   JOURNAL = {Commun. Stoch. Anal.},
  FJOURNAL = {Communications on Stochastic Analysis},
    VOLUME = {7},
      YEAR = {2013},
    NUMBER = {2},
     PAGES = {203--225},
   MRCLASS = {91G10 (60G44 60H30 62P05 91B30)},
  MRNUMBER = {3092230},
MRREVIEWER = {Yong Hyun Shin},
       DOI = {10.31390/cosa.7.2.03},
       URL = {https://doi.org/10.31390/cosa.7.2.03},
}

@article {peng90,
    AUTHOR = {Peng, Shi Ge},
     TITLE = {A general stochastic maximum principle for optimal control
              problems},
   JOURNAL = {SIAM J. Control Optim.},
  FJOURNAL = {SIAM Journal on Control and Optimization},
    VOLUME = {28},
      YEAR = {1990},
    NUMBER = {4},
     PAGES = {966--979},
      ISSN = {0363-0129},
   MRCLASS = {49K45 (93E20)},
  MRNUMBER = {1051633},
MRREVIEWER = {Henryk G\'{o}recki},
       DOI = {10.1137/0328054},
       URL = {https://doi.org/10.1137/0328054},
}

@article {peng93,
    AUTHOR = {Peng, Shige},
     TITLE = {Backward stochastic differential equations and applications to
              optimal control},
   JOURNAL = {Appl. Math. Optim.},
  FJOURNAL = {Applied Mathematics and Optimization},
    VOLUME = {27},
      YEAR = {1993},
    NUMBER = {2},
     PAGES = {125--144},
      ISSN = {0095-4616},
   MRCLASS = {49K45 (60H10 93E20)},
  MRNUMBER = {1202528},
MRREVIEWER = {N. U. Ahmed},
       DOI = {10.1007/BF01195978},
       URL = {https://doi.org/10.1007/BF01195978},
}

@article {pennerreveillac15,
    AUTHOR = {Penner, Irina and R\'{e}veillac, Anthony},
     TITLE = {Risk measures for processes and {BSDE}s},
   JOURNAL = {Finance Stoch.},
  FJOURNAL = {Finance and Stochastics},
    VOLUME = {19},
      YEAR = {2015},
    NUMBER = {1},
     PAGES = {23--66},
      ISSN = {0949-2984},
   MRCLASS = {91B30 (60G07 60G40 60H10 91B16)},
  MRNUMBER = {3292124},
MRREVIEWER = {Lingjiong Zhu},
       DOI = {10.1007/s00780-014-0243-x},
       URL = {https://doi.org/10.1007/s00780-014-0243-x},
}

@article {pirvu07,
    AUTHOR = {Pirvu, Traian A.},
     TITLE = {Portfolio optimization under the value\:\!-\:\!at\:\!-risk constraint},
   JOURNAL = {Quant. Finance},
  FJOURNAL = {Quantitative Finance},
    VOLUME = {7},
      YEAR = {2007},
    NUMBER = {2},
     PAGES = {125--136},
      ISSN = {1469-7688},
   MRCLASS = {91B28 (60H10 60H30 90C39 93E20)},
  MRNUMBER = {2325660},
MRREVIEWER = {Guohe Deng},
       DOI = {10.1080/14697680701213868},
       URL = {https://doi.org/10.1080/14697680701213868},
}

@article {pollak68,
    AUTHOR = {Pollak, R. A.},
     TITLE = {Consistent planning},
   JOURNAL = {Rev. Econ. Stud},
  FJOURNAL = {},
    VOLUME = {35},
      YEAR = {1968},
    NUMBER = {},
     PAGES = {185--199},
      ISSN = {},
   MRCLASS = {},
  MRNUMBER = {},
MRREVIEWER = {},
       DOI = {},
       URL = {https://www.jstor.org/stable/2296548},
}

@article {strotz55,
    AUTHOR = {Strotz, R. H.},
     TITLE = {Myopia and inconsistency in dynamic utility maximization},
   JOURNAL = {Rev. Econ. Stud.},
  FJOURNAL = {},
    VOLUME = {23},
      YEAR = {1955},
    NUMBER = {},
     PAGES = {165--180},
      ISSN = {},
   MRCLASS = {},
  MRNUMBER = {},
MRREVIEWER = {},
       DOI = {},
       URL = {https://www.jstor.org/stable/2295722},
}

@article {weiyongyu17,
    AUTHOR = {Wei, Qingmeng and Yong, Jiongmin and Yu, Zhiyong},
     TITLE = {Time\:\!-inconsistent recursive stochastic optimal control
              problems},
   JOURNAL = {SIAM J. Control Optim.},
  FJOURNAL = {SIAM Journal on Control and Optimization},
    VOLUME = {55},
      YEAR = {2017},
    NUMBER = {6},
     PAGES = {4156--4201},
      ISSN = {0363-0129},
   MRCLASS = {93E20 (49N70 91A23 91A65)},
  MRNUMBER = {3738841},
MRREVIEWER = {Annalisa Cesaroni},
       DOI = {10.1137/16M1079415},
       URL = {https://doi.org/10.1137/16M1079415},
}

@Article{yong99,
  title =	 {Infinite dimensional optimal control theory},
  author =	 {Yong, J.},
  journal =	 {IFAC Proceedings Volumes, Elsevier},
  year =	 {1999}
}

@book {yongzhou99,
    AUTHOR = {Yong, Jiongmin and Zhou, Xun Yu},
     TITLE = {Stochastic controls},
    SERIES = {Applications of Mathematics (New York)},
    VOLUME = {43},
      NOTE = {Hamiltonian systems and HJB equations},
 PUBLISHER = {Springer-Verlag, New York},
      YEAR = {1999},
     PAGES = {xxii+438},
      ISBN = {0-387-98723-1},
   MRCLASS = {93-02 (49K45 49L20)},
  MRNUMBER = {1696772},
MRREVIEWER = {Tamer Ba\c{s}ar},
       DOI = {10.1007/978-1-4612-1466-3},
       URL = {https://doi.org/10.1007/978-1-4612-1466-3},
}

@article {yong12,
    AUTHOR = {Yong, Jiongmin},
     TITLE = {Time\:\!-inconsistent optimal control problems and the equilibrium
              {HJB} equation},
   JOURNAL = {Math. Control Relat. Fields},
  FJOURNAL = {Mathematical Control and Related Fields},
    VOLUME = {2},
      YEAR = {2012},
    NUMBER = {3},
     PAGES = {271--329},
      ISSN = {2156-8472},
   MRCLASS = {49K45 (35R60 49L25 49N70 91B50 93E20)},
  MRNUMBER = {2991570},
MRREVIEWER = {Jamsheed Shorish},
       DOI = {10.3934/mcrf.2012.2.271},
       URL = {https://doi.org/10.3934/mcrf.2012.2.271},
}

@inproceedings {yong14,
    AUTHOR = {Yong, Jiongmin},
     TITLE = {Time-inconsistent optimal control problems},
 BOOKTITLE = {Proceedings of the {I}nternational {C}ongress of
              {M}athematicians---{S}eoul 2014. {V}ol. {IV}},
     PAGES = {947--969},
 PUBLISHER = {Kyung Moon Sa, Seoul},
      YEAR = {2014},
   MRCLASS = {93E20 (49K45 93-02)},
  MRNUMBER = {3751160},
MRREVIEWER = {Francisco Jos\'{e} Silva \'{A}lvarez},
}

@article {zhuo18,
    AUTHOR = {Zhuo, Yu},
     TITLE = {Maximum principle of optimal stochastic control with terminal
              state constraint and its application in finance},
   JOURNAL = {J. Syst. Sci. Complex.},
  FJOURNAL = {Journal of Systems Science \& Complexity},
    VOLUME = {31},
      YEAR = {2018},
    NUMBER = {4},
     PAGES = {907--926},
      ISSN = {1009-6124},
   MRCLASS = {93E20 (49K45 91G10)},
  MRNUMBER = {3788948},
MRREVIEWER = {Mokhtar Hafayed},
       DOI = {10.1007/s11424-018-6212-2},
       URL = {https://doi.org/10.1007/s11424-018-6212-2},
}

@article{Ha21MCRF,
  title={Extended backward stochastic Volterra integral equations and their applications to time-Inconsistent stochastic recursive control problems},
  author={Hamaguchi, Yushi},
  journal={Mathematical Control \& Related Fields},
  volume={11},
  number={2},
  pages={433},
  year={2021},
  publisher={American Institute of Mathematical Sciences}
}

@Misc{amsmath,
  author =	 {{American Mathematical Society}},
  title =	 {User's Guide for the \texttt{amsmath} Package
                  (Version 2.0)},
  url =		 {ftp://ftp.ams.org/pub/tex/doc/amsmath/amsldoc.pdf},
  urldate =	 {2015-07-30},
  year =	 2002}





\end{document}